\documentclass{sn-jnl}

\usepackage{graphicx}
\usepackage{multirow}
\usepackage{amsmath,amssymb,amsfonts}
\usepackage{amsthm}
\usepackage{mathrsfs}
\usepackage[title]{appendix}
\usepackage{xcolor}
\usepackage{textcomp}
\usepackage{manyfoot}
\usepackage{booktabs}
\usepackage{algorithm}
\usepackage{algorithmicx}
\usepackage{algpseudocode}
\usepackage{listings}
\usepackage{xargs}   
\usepackage{enumerate}
\usepackage{tikz-cd}
\usepackage{mathtools}
\usepackage{breqn}
\usepackage{array}
\usepackage{booktabs}
\usepackage{tabularx}
\usepackage{makecell}
\usepackage{longtable}
\usepackage{multirow}

\theoremstyle{thmstyleone}%
\newtheorem{theorem}{Theorem}
\newtheorem{proposition}[theorem]{Proposition}%
\newtheorem{lemma}[theorem]{Lemma}%
\theoremstyle{thmstyletwo}%
\newtheorem{example}{Example}%
\newtheorem{remark}{Remark}%
\theoremstyle{thmstylethree}%
\newtheorem{definition}{Definition}%

\newcommand{\set}{\mathbf{Set}_{k}}
\newcommand{\setnd}{\mathbf{Set}_{k}^{ND}}
\newcommand{\mc}{\mathbf{Met}_{k}^{C}}
\newcommand{\mcnd}{\mathbf{Met}_{k}^{CND}}
\newcommand{\mcndk}[1]{\mathbf{Met}_{#1}^{CND}}
\newcommand{\setndk}[1]{\mathbf{Set}_{#1}^{ND}}

\newcommand{\M}{\Sigma}
\newcommand{\F}{\M\otimes}
\newcommand{\Fw}{\M\otimes_{\w}}
\newcommand{\Fww}[1]{\M\otimes_{#1}}
\newcommand{\Fwi}[1]{\M^{#1}\otimes_{\w}}
\newcommand{\Fi}[1]{\M^{#1}\otimes }
\newcommand{\Mw}{\M\times_{\w}}

\newcommand{\A}{\left(A, a_1 ,\cdots , a_k\right)}
\newcommand{\X}{\left(X, x_1 ,\cdots , x_k\right)}
\newcommand{\B}{\left(B, b_1, \cdots, b_k \right)}
\newcommand{\G}{\mathbb{G}}
\newcommand{\g}{\tilde{g}^{-1}}

\newcommand{\dm}[1]{d_{\M \times_{\w} #1}}
\newcommand{\dqm}[1]{d_{\M \otimes_{\w} #1}}
\newcommand{\f}{\phi}
\newcommand{\tf}{\tilde{\Phi}}
\newcommand{\w}{\rho}
\newcommand{\s}{\mu}
\newcommand{\Ss}{\mathbb{S}}
\newcommand{\gPoints}{\pi}
\newcommand{\ft}{\widetilde f}

\begin{document}

\title[Presenting Point-Glued Spaces as Final Coalgebras via Metric Completion]{Presenting Point-Glued Self-Similar Objects as Final Coalgebras in Metric Spaces via Completion: A Categorical Hutchinson Theorem}

\author*[1]{\fnm{Jayampathy} \sur{Ratnayake}}\email{jratnayake@sjp.ac.lk}

\author[2]{\fnm{Annanthakrishna} \sur{Manokaran}}\email{annanth@univ.jfn.ac.lk}

\affil*[1]{\orgdiv{Department of Mathematics}, \orgname{University of Sri Jayewardenepura}, \orgaddress{\city{Nugegoda}, \country{Sri Lanka}}}

\affil[2]{\orgdiv{Department of Mathematics and Statistics}, \orgname{University of Jaffna}, \orgaddress{\city{Jaffna}, \country{Sri Lanka}}}

\abstract{
We study self-similar spaces obtained by repeatedly gluing finitely many scaled copies of a pointed metric space along prescribed boundary points. A finite gluing specification defines endofunctors on suitable categories of pointed sets and pointed metric spaces. Under non-degeneracy and metric admissibility assumptions, the initial algebra of the set functor carries a canonical metric. Its Cauchy completion is a final coalgebra for the metric functor and, after forgetting the metric, also for the underlying set functor.

We compare this construction with Iterated Function Systems. When the orbit algebra is injective and its comparison with the initial algebra is bi-Lipschitz, the attractor is isomorphic to the completed final coalgebra. This gives a categorical analogue of Hutchinson's theorem for overlaps captured by finitely many distinguished points. Examples include the unit interval, Cantor spaces, the Sierpiński gasket and $p$-adic integers.
}



\maketitle

\section{Introduction}\label{introduction}
Self-similar spaces such as the unit interval and the Sierpiński gasket can be constructed by gluing finitely many scaled copies of themselves, with overlaps prescribed at distinguished points. These spaces also arise as attractors of iterated function systems (IFSs). For example, the unit interval satisfies
\[
    [0,1]=\f_1([0,1])\cup\f_2([0,1]),  \qquad  \f_1(x)=\frac{x}{2}, \quad  \f_2(x)=\frac{x+1}{2}.
\]
The two images are copies of $[0,1]$ scaled by $1/2$ and joined end to end at $\f_1(1)=\f_2(0)=1/2$.

Hutchinson's theorem \cite{Hut} associates a unique nonempty compact attractor with a finite family of contractions on a complete metric space. The attractor is the fixed point of the induced hyperspace operator and can also be recovered
concretely as the closure of the ``finitely specified" orbit: the images of the branch fixed points under all finite compositions of the contractions. For the interval IFS above, this orbit consists of dyadic rationals in $[0,1]$.

A main point of the present paper is that this picture has a categorical analog.  In our setting, the role of the fixed points and their finite words is played by the initial algebra of a point gluing endofunctor.  The role of the closure is played by Cauchy completion. The resulting complete object is then described as the largest fixed point of the functor: the final coalgebra characterized by its universal property. Thus, our result may be regarded as Hutchinson-type existence and uniqueness theorems in a categorical setting: the limiting self-similar space is obtained as a fixed point by completing the finite-stage object and is universal with respect to the associated coalgebra structure. We summarise this analogy below, in Table \ref{tab:comparisonHutCat}.
\begin{table}[b]
    \centering
\renewcommand{\arraystretch}{1.25}
\begin{tabular}{>{\raggedright\arraybackslash}p{0.43\textwidth}|>{\raggedright\arraybackslash}p{0.43\textwidth}}
    \textbf{Hutchinson/IFS picture} & \textbf{Categorical gluing picture} \\
    \hline
    Finite family of contractions & Finite pointed gluing data, a contraction factor and an associated endofunctor \\
    Finite compositions of contractions & Iterates of the gluing functor \\
    Images of the branch fixed points under finite compositions & The initial algebra \\
    Closure of these finite approximants & Cauchy completion of the initial algebra \\
    Unique compact attractor & Final coalgebra of the gluing functor \\
    Fixed-point/self-similarity equation & Structure isomorphism $\mathbb{G}\cong \Fw \mathbb{G}$\\
    \botrule
\end{tabular}
\caption{Hutchinson and categorical comparison}
\label{tab:comparisonHutCat}
\end{table}

This work builds on categorical treatments of self-similarity, including Freyd's characterization of the unit interval as a final coalgebra \cite{freyd:real}, and Leinster's general framework of gluing equations   \cite{Lein}. Leinster describes self-similarity using equational systems, specified by finite non-degenerate modules on small categories, and their universal solutions \cite[Definitions~2.8--2.9]{Lein}. His IFS realization   \cite[Example~10.11]{Lein} is a direct antecedent of our IFS comparison theorem. 

Using Leinster’s framework, Hasuo, Jacobs and Niqui \cite{coalge_representation} develop symbolic representations of fractals by equivalence classes of infinite streams, with the identifications accounting for gluing. Under an injectivity hypothesis, they prove that the canonical denotation map is injective, so each point in its image has a unique symbolic representation modulo the prescribed identifications. This establishes a bijection between the symbolic final coalgebra and the metric attractor. The present paper develops an initial-algebra completion construction for finite pointed gluing specifications, with explicit metric admissibility and quantitative conditions for comparison with IFS attractors.

In a complementary direction, Schmid, Noquez and Moss \cite{SchmidNoquezMoss2025} relate IFS attractors to process calculi, interpreting process expressions using contractions on complete metric spaces and providing a sound and complete axiomatization of fractal equivalence.

The passage from an initial algebra to a final coalgebra by completion has a broader categorical history. An order-theoretic approach was developed by Adámek \cite{Adamek2002Ideal}, who realized final coalgebras as ideal completions of suitably ordered initial algebras under appropriate hypotheses on the functor. For finitary set functors $F$ with $F\varnothing\neq\varnothing$ that preserve limits of inverse $\omega$-chains, Barr \cite{Bar} identified the terminal coalgebra with the Cauchy completion of the initial algebra in a canonical ultrametric. Without the limit-preservation assumption, Ad\'amek \cite{Adamek2019Completion,Adamek2021Completion} proved that the two objects carry canonical ultrametrics with a common Cauchy completion; the latter extends this result to precontinuous set functors.  A systematic account of these results appears in \cite[Chapter~13]{AMM}.

The present work builds on earlier metric presentations of the unit interval \cite{Bhat,BhatTACL}, the Sierpiński gasket \cite{Sierpinski_gasket}, the Sierpiński carpet \cite{Sierpinsky_carpet}, and the $p$-adic integers \cite{Mano}. Our construction concerns pointed objects equipped with quotient metrics induced by scaled finite gluing. Under non-degeneracy and metric admissibility, the set-theoretic initial algebra carries a canonical metric whose Cauchy completion is a final coalgebra for the  set functor and its continuous metric counterpart. These quotient metrics need not be ultrametric, allowing connected examples such as the unit interval. We also establish quantitative conditions under which this completion can be identified with an IFS attractor.
\begin{equation} \nonumber
    \underbrace{(G,g)}_{\text{set-theoretic initial algebra}}
    \quad\xrightarrow{\;\substack{\text{canonical metric}\\ \text{and} \;\\{\;\;\text{completion}\;}}}\quad
    \underbrace{(\G,\g)}_{
        \substack{
            \text{final in the set setting}\\
            \text{and, with the canonical metric,}\\
            \text{final in the metric setting}
        }
    }
\end{equation}

\section{Results}
We start with a finite gluing specification $\Sigma=(M,D,R)$.  Here $M$ indexes the copies used in each stage, $D$ records the output distinguished points, and the relation $R$ prescribes boundary identifications. Appendix~\ref{Subsec:RigorousModuleConstruction} identifies the corresponding gluing functor with the tensor functor of an associated finite module (as in \cite{Lein}) on the two-object category presenting $k$-pointed sets. With a scaling factor $ 0< \rho <1$,  it defines endofunctors $\F-$ and $\Fw-$ on suitable categories of non-degenerate pointed sets ($\setnd$) and metric spaces ($\mcnd$), respectively. Applied to a pointed metric space $X$, this functor forms finitely many $\rho$-scaled copies of $X$ and identifies distinguished points according to the prescribed gluing relation.  Non-degeneracy keeps the output distinguished points distinct, while metric admissibility (cf. Section \ref{subsec_MetricAdmissibilityAndCopyseparation}) allows the initial chain to be realized by isometric embeddings. The colimit of the initial chain yields the initial algebra $G$ in $\setnd$ and we give it the metric inherited from the isometric embeddings of the chain. Write $\G=C(G)$ for its Cauchy completion. Our first main theorem shows that this completion carries a natural coalgebra structure and is the final coalgebra of the corresponding set gluing functor. 

\begin{theorem}[Final coalgebra in $\setnd$]
\label{Thm_FinalCoalgebraInSet}
     Assume that the finite gluing specification $\Sigma=(M,D,R)$ is non-degenerate and  $(\Sigma,\w)$ is metrically admissible for some $\w \in (0,1)$. 
     Then the initial algebra $(G, g)$ of the corresponding set functor $\F-$ on $\setnd$ is the colimit of the initial chain associated with the initial object of $\setnd$, and its completion $(\G, \g)$ with respect to the canonical metric is the final coalgebra of $\F-$ on $\setnd$.
\end{theorem} 

Finality in the continuous metric category follows by proving that the unique set-theoretic mediating map is continuous whenever the coalgebra structure is continuous.

\begin{theorem}[Final coalgebra in $\mcnd$]
\label{Thm_FinalCoalgebraInMet}
    Assume that $(\Sigma,\w)$ is metrically admissible. Let $(X,e)$ be an $\Fw-$coalgebra on $\mcnd$; that is, $e:X\rightarrow \Fw X$ is continuous.  Then, the unique map $f:X\rightarrow \G$ obtained in Theorem \ref{Thm_FinalCoalgebraInSet} is continuous.  Consequently, $(\G,\g)$ is the final coalgebra of $\Fw-$ on $\mcnd$.
\end{theorem}

Finally, for a hyperbolic IFS, we prove a categorical comparison with Hutchinson's theorem. The associated finite point-gluing specification $\Sigma$ (Definition~\ref{Def_RelationForIFS}) gives rise to our gluing functor, and the orbit of the branch fixed points carries a natural algebra structure (Proposition~\ref{prop_OrbithAlgebraAndFiniteStages}). Under the hypotheses below, the classical attractor can then be given a coalgebra structure, which can be identified with the final coalgebra of Theorem \ref{Thm_FinalCoalgebraInMet}. This provides a categorical refinement of Hutchinson's attractor construction under the stated compatibility conditions.

\begin{theorem}[Categorical Hutchinson comparison]
\label{Thm_F_isTheFinalCoalgebra}
Let $(\f_1,\ldots,\f_k)$, with $k>1$, be a finite IFS on a nonempty complete one-bounded metric space, with a common Lipschitz bound $0<\lambda<1$ and pairwise distinct fixed points $\alpha_1,\ldots,\alpha_k$. Let $O$ be the orbit of the fixed points under finite compositions of the branches and write $\mathcal F=\overline O$ for the attractor.

The associated gluing specification $\Sigma=(M,D,R)$, with $M=\{1,\ldots,k\}$, defined in Definition~\ref{Def_RelationForIFS}, is non-degenerate.

Suppose that $0<\rho<1$ is such that $(\Sigma,\rho)$ is metrically admissible. Let $(G,g)$ be the initial algebra of the associated set functor, equipped with its canonical metric. Assume that the orbit algebra $\f:\Fw O\longrightarrow O$, given by $\f(i\otimes x)=\f_i(x)$, is injective and that the canonical
algebra morphism $f:G\longrightarrow O$ is bi-Lipschitz. Then:
    \begin{enumerate}[(i)]
        \item $(O,\f)$ is isomorphic to the algebra $(G,g)$ on $\mcnd$.
       \item Transport along $\kappa_O$ of the completed orbit algebra,
        \begin{equation} 
        \nonumber
            \widetilde{\f} : \Fw\mathcal F\longrightarrow\mathcal F,
        \end{equation}
        is an isomorphism, and its inverse
        \begin{equation} 
        \nonumber
            {\widetilde{\f}}^{-1}: \mathcal F\longrightarrow\Fw\mathcal F
        \end{equation}
        is a final coalgebra. Equivalently, $(\mathcal F,\widetilde{\f}^{-1}) \cong (\G,\g)$.
        \item  After forgetting the metric structure, $(\mathcal F,\widetilde{\f}^{-1})$ is also final for the underlying set functor on $\setnd$.
    \end{enumerate}
Here, $\kappa_O:C(O)\longrightarrow \overline O=\mathcal F$ is the canonical isometric isomorphism given by $\kappa_O([x_n])=\lim_{n\to\infty}x_n$, and $(\G,\g)$ is the final coalgebra obtained by completing the set theoretic initial algebra $(G,g)$ with respect to its canonical metric.
\end{theorem}

The paper is organized as follows.  Section~\ref{Sec:Definitions} introduces the pointed categories and gluing endofunctors. Sections~\ref{Sec_InitialAlgebraOfF} and~\ref{Sec:Completion} construct the initial algebra, its canonical metric, and its completion. Section~\ref{Sec:FinalCoalgebra} proves finality, and Section~\ref{Sec:IFS} develops the IFS comparison. Section~\ref{Sec:Examples} contains examples, including Cantor-type spaces, the Sierpinski gasket, $p$-adic integers, and nested radicals of $2$. The carpet example (Example~\ref{Ex_Sierpinski_Carpet_IFS}) illustrates the limitations of point gluing when geometric overlaps extend beyond the prescribed boundary points. We take the unit interval as a running example to demonstrate the constructions.

\section{Basic Definitions and Constructions}
\label{Sec:Definitions}

For a positive integer $n$, write $[n]=\{1,2,\cdots, n\}$.

\subsection{Categories of $k$-pointed sets and metric spaces}

\subparagraph{Categories of $k$-pointed sets} A $k$-pointed set, where $k\in\mathbb{N}$, is a set $A$ together with $k$ distinguished elements $\left( a_1,\cdots, a_k \right)\in A^k$. Such a set will be denoted as $\left( A\, ,\,a_1,\cdots,a_k\right)$. We let $A_\partial=\{a_j \,\,:\,\, j \in [k] \}$, the boundary of $A$, to be the set of distinguished points. Let $\set$ be the category whose objects are $k$-pointed sets and morphisms are pointed functions. Thus a morphism $f: \A \to \B$ in $\set$ is a function satisfying $f(a_j)=b_j$ for all $j \in [k]$. 

When it is clear from context, we often omit the distinguished points from the notation and write only the underlying set. The initial $k$-pointed set is
\begin{equation}
\nonumber
    I=([k],1,\ldots,k).
\end{equation}
For $A=(A,a_1,\ldots,a_k)$, write $!_A:I\to A$ for the unique pointed map, given by $!_A(j)=a_j$. Metric categories are defined similarly.

\subparagraph{Categories of $k$-pointed metric spaces}   A $k$-pointed metric space is a $k$-pointed set $A$ together with a one-bounded metric $d$; i.e., $d(a, a')\leq 1$ for all $a, a' \in A$. Let $\mc$ be the category of $k$-pointed metric spaces whose morphisms are continuous functions that preserve the distinguished elements. 

\subparagraph{Non-degenerate subcategories} A $k$-pointed set $\A$ is non-degenerate if all the distinguished points are pairwise distinct. Let $\set^{ND}$ be the full subcategory of $\set$, whose objects are the non-degenerate $k$-pointed sets. Non-degeneracy in metric setting means that the underlying $k$-pointed set is non-degenerate. Write $\mcnd$ for the full subcategory of non-degenerate objects in $\mc$.

For a $k$-pointed metric space $\A$, we let
    \begin{equation} 
    \nonumber
        \Delta(A)=
        \begin{cases}
            1 & k=1,\\[4pt]
            \displaystyle \min\left\{d(a_{j}, a_{j'}) \,\,:\,\, 1\leq j<j'\leq k\right\} & k\geq 2.
        \end{cases}
    \end{equation}

For $k=1$, the convention $\Delta(A)=1$ reflects the facts that every one-pointed metric space is non-degenerate and all metrics considered are one-bounded. Note that a $k$-pointed metric space $A$ is non-degenerate if and only if $\Delta(A)>0$.

\subsection{Construction of the Endofunctors}
Leinster \cite{Lein}  describes self-similar spaces obtained by gluing copies of spaces as prescribed by a pair $(\mathcal A,\mathcal M)$ consisting of a small category $\mathcal A$  and a finite non-degenerate module  $\mathcal M$. Here, the module encodes the copying and gluing data. In \cite{coalge_representation} authors use the term \emph{combinatorial specification} for such a pair. 

In the finite pointed setting, we encode the copies, output distinguished points, and boundary identifications by a finite gluing specification $\Sigma=(M,D,R)$. This description makes the copying and gluing data explicit. Appendix~\ref{Subsec:RigorousModuleConstruction} constructs the associated finite module and identifies its tensor functor with the concrete gluing functor. 

\begin{definition}[Finite gluing specification over points]
\label{Def_FiniteGluingSpecification}
    Fix $k\in\mathbb N$.  A \emph{finite $k$-pointed gluing specification} (or simply a \emph{gluing specification}) is a triple
    \begin{equation} 
    \nonumber
        \Sigma=(M,D,R),
    \end{equation}
    where
    \begin{enumerate}
        \item $M=[\s]$, where $\s\geq1$, is the ordered set of copy indices used in one step of the construction;
        \item $D=\left((i_1,j_1),\cdots,(i_k,j_k)\right)$ is a $k$-tuple of \underline{distinct} elements of $M\times I$; and
        \item $R$ is a reflexive and symmetric binary relation on $M\times I$.
    \end{enumerate}
    Write $\overline R$ for the transitive closure, equivalently the equivalence closure, of $R$. We also use $D$ for its underlying set. Write $d_r=(i_r, j_r)$ for the $r$-th \textit{output label}, i.e., the $r$-th entry of $D$.
\end{definition}

Fix a finite gluing specification $\Sigma=(M,D,R)$. For a $k$-pointed set $\A$, write
\begin{equation} 
\nonumber
    \M\times A=(M\times A,(i_1,a_{j_1}),\ldots,(i_k,a_{j_k})).
\end{equation}
This forms $\mu$ copies of $A$, with the distinguished points selected as prescribed by $D$. Specifically, if $d_r=(i_r, j_r)$ is the $r$-th element of $D$, then the $r$-th distinguished point of $\M \times A$ is $(i_r,a_{j_r})$, the $j_r$-th distinguished point in the ${i_r}$-th copy of $A$. For a pointed map $f:A\to B$, put $(M\times f)(m,a)=(m,f(a))$. These definitions make  $\M\times -$ an endofunctor on $\set$.

In particular, the specification determines the non-degenerate $k$-pointed set
\begin{equation}
\nonumber
    \mathcal P=\M\times I=(M\times I,d_1,\ldots,d_k).
\end{equation}
The unique pointed map $!_A:I\to A$ induces the pointed map
\begin{equation}
\nonumber
    \gPoints_A=\M\times !_A:\mathcal P\longrightarrow\M\times A,
    \qquad \gPoints_A(i,j)=(i,a_j).
\end{equation}

The map $\gPoints_A$ determines a quotient of $\M \times A$ on $\set$ as follows.

\begin{definition}
\label{Def_gluingRel}
    Given a gluing specification $\Sigma=(M,D,R)$, with $k=|D|$ and $\mu=|M|$, and a $k$-pointed set $\left( A, a_1, \cdots, a_k \right)$, $R$ induces a relation $\sim_A$ on $M\times A$ given by 
    \begin{equation} \nonumber
       \gPoints_A(i,j)\sim_A\gPoints_A(i',j')
       \quad\text{whenever}\quad
       (i,j)R(i',j').
    \end{equation}
    That is, $\sim_A\, =\, \{\left(\left(i, a_j \right)\,,\,\left(i', a_{j'} \right)\right)\,:\, \,\,\,\left(\left(i,j\right)\,,\, \left(i', j'\right)\right)\in R \}$.
    
    Let $\approx_A$ be the equivalence closure of $\sim_A$, and 
    \begin{equation} 
    \nonumber
        M \otimes A = M \times A/\approx_A
    \end{equation} 
    be the set of equivalence classes with respect to $\approx_A$.      We denote the equivalence class of an element $(m, a)\in M \times A$ by $m\otimes a \in M \otimes A$. The set $M\otimes A$ is lifted to a $k$-pointed set by choosing $i_r\otimes a_{j_r}$ as the $r$-th distinguished point, where $d_r=(i_r, j_r)\in D$ and $r \in [k]$. We denote this $k$-pointed set by $\F A$. I.e.,
    \begin{equation} 
    \nonumber
        \F A = \Big(M \otimes A, i_1\otimes a_{j_1}, \cdots, i_k\otimes a_{j_k} \Big).
    \end{equation}
\end{definition}

Put $\mathcal C=\F I$. Since $\gPoints_I$ is the identity on $M\times I$, we have $\approx_I=\overline R$, and hence
\begin{equation}
\nonumber
    \mathcal C=\bigl((M\times I)/\overline R,C_1,\ldots,C_k\bigr),
    \qquad C_r=i_r\otimes j_r=[d_r]_{\overline R}.
\end{equation}
The quotient map $\mathcal P\to\mathcal C$ is pointed. The unique pointed map from $I$ to $\mathcal C$ is
\begin{equation}
\nonumber
    \eta:I\longrightarrow\mathcal C=\F I,\qquad \eta(r)=C_r.
\end{equation}
We also use $\mathcal P$ and $\mathcal C$ for their underlying sets.

Intuitively, the construction of $M\otimes A$ is understood as gluing $|M|=\s$ copies of $A$ in $\M \times A$, where gluing takes place only along the distinguished points in various copies of $A$ in $M\times A$ and governed by the relation $R$. Only the equivalence closure of the prescribed boundary pairs affects the quotient. Allowing $R$ to be reflexive and symmetric, rather than requiring transitivity, lets examples specify generating pairs. 

We give the unit interval as a running example to demonstrate the key construction.
\begin{example}[Unit interval: gluing data]
\label{Ex_unitInterval}
    Take $k=2$, $M=[2]$, and $D=((1,1),(2,2))$. Let $R$ be the equivalence relation generated by $(1,2)\mathrel{R}(2,1)$. 
    
    The specification $\Sigma = (M, D, R)$ prescribes the following gluing of a bipointed set $(A,a_1,a_2)$. Take two copies $\{1\} \times A$ and $\{2\} \times A$ of $A$ and join them end-to-end by the identification $(1,a_2)\sim_A (2,a_1)$. The quotient $\F A$ is given the distinguished points $1\otimes a_1$ and $2\otimes a_2$ (see \cite{freyd:real,Bhat}).
    \begin{equation}
    \nonumber
        \F A = \Bigl( \bigl( [2] \times A \bigr)/\langle(1,a_2)\sim_A (2,a_1)\rangle\;,\; 1\otimes a_1\;,\;2\otimes a_2 \Bigr)
    \end{equation}
\end{example}
    
The construction of Definition \ref{Def_gluingRel} has the following  properties.

\begin{proposition}[Elementary properties of the gluing construction]
\label{Prop_ElementryPropertiesOfGluingConstruction}
For a finite gluing specification $\Sigma=(M,D,R)$:
\begin{enumerate}[(i)]
    \item  \label{Lem_EndoOnSet_k}
        The functor $\M\times-$ induces an endofunctor $\F-$ on $\set$ given by the assignment
        \begin{equation} \nonumber
            (\F f)(m\otimes a)=m\otimes f(a).
        \end{equation}
    \item \label{Lem:EquivalenceClosureInvariance} 
        Replacing $R$ by $\overline R$ leaves the pointed quotient         $M\otimes A$ unchanged.
    \item \label{Lem_InjectionsReflectGluing}
         If $f:A\to B$ is injective, then
        \begin{equation} \nonumber
            (M\times f)(z)\approx_B(M\times f)(z') \quad \Longleftrightarrow \quad z\approx_A z', \quad  \text{ for all }\; z,z'\in M\times A.
        \end{equation}
    \item \label{Lem_oneone}
        Consequently, if $f: A \to B$ is injective, then so is $\F f$.
\end{enumerate}
\end{proposition}

\begin{proof}
Pointed maps preserve the generating gluing pairs, so $M\times f$ induces the pointed map $\F f$ on the quotients. The functor laws follow directly from the displayed formula, proving (i).  

For (ii), we note that an $\overline R$-pair is connected by an $R$-zig-zag, so its image already belongs to $\approx_A$. 

Every endpoint of a nontrivial gluing step lies in $M\times B_\partial$, and hence in the image of $M\times f$. Injectivity therefore allows each gluing zig-zag between image points to be pulled back uniquely. Thus (iii) holds.

For (iv), equality of quotient classes is thus reflected, which proves injectivity of $\F f$.
\end{proof}

Applying $\F-$ to $!_A:I\to A$ gives the pointed map
\begin{equation}
\nonumber
    \overline{\gPoints}_A=\F(!_A):\mathcal C\longrightarrow\F A,
    \qquad \overline{\gPoints}_A(i\otimes j)=i\otimes a_j.
\end{equation}
Its image is the finite set of quotient classes represented by boundary points:
\begin{equation}
\label{eqn:BoundaryImage}
    \mathcal C_A=\operatorname{im}(\overline{\gPoints}_A)
    =\{i\otimes a_j:i\in M,\ j\in I\}\subseteq\F A.
\end{equation}
If $A$ is non-degenerate, then $!_A$ is injective, so Proposition~\ref{Prop_ElementryPropertiesOfGluingConstruction}-(\ref{Lem_oneone}) makes $\overline{\gPoints}_A$ injective. Thus $\mathcal C$ is canonically identified, as a pointed set, with $\mathcal C_A$.

Note that $\sim_A$ is ``essentially" a relation on $M\times A_\partial$. Specifically, the only relations other than reflexive ones in the equivalence relation $\approx_A$ are of the form $(m, a_i)\approx_A (m', a_{i'})$, where $a_i, a_{i'}\in A_\partial$ (cf. Lemma \ref{Lem:BoundaryClassCriterion} and \ref{Lem:FiniteQuotientFibres}). Thus, only the equivalence relation generated by the prescribed gluing pairs affects the quotient. 

If distinguished points of $A$ coincide, different boundary labels may denote the same point. For non-degenerate $A$, the boundary labels faithfully record the prescribed gluing classes, as follows.

\begin{lemma}[Boundary-class criterion]
\label{Lem:BoundaryClassCriterion}
    Let $ A$ be a $k$-pointed set.
    \begin{enumerate}[(i)]
        \item For all $p,q\in \mathcal P$,
        \begin{equation} 
        \nonumber
            p\mathrel{\overline R}q \quad \Longrightarrow \quad       \gPoints_A(p) \mathrel{\approx_A} \gPoints_A(q).
        \end{equation}
        \item The converse holds if $ A$ is non-degenerate: for all $p,q \in \mathcal P$
        \begin{equation} 
        \nonumber
           \gPoints_A(p) \mathrel{\approx_A} \gPoints_A(q) \quad \Longleftrightarrow \quad p\mathrel{\overline R}q.
       \end{equation}
    \end{enumerate}
\end{lemma}

\begin{proof}
For $p,q\in\mathcal P$, the relation $p\mathrel{\overline R}q$ is equivalent to equality of their classes in $\mathcal C=\F I$. Applying $\overline{\gPoints}_A$ gives (i). If $A$ is non-degenerate, this map is injective, so equality of the image classes implies equality in $\mathcal C$. This proves (ii).
\end{proof}

For $C\in\mathcal C$, define the finite set of boundary representatives
\begin{equation} 
\nonumber 
   B_C=\{(i,a_j):(i,j)\in C\}\subseteq M \times A_\partial,
\end{equation}
and for $u\in M \otimes A$, put
\begin{equation} 
\nonumber
   \operatorname{Rep}(u)= \{ (m,a) \in \M\times A : m\otimes a = u\}.
\end{equation} 

\begin{lemma}[Finite quotient fibres]
\label{Lem:FiniteQuotientFibres}
    For any $u\in \F A$, $\operatorname{Rep}(u)$ is finite.  More precisely, for a $k$-pointed set $A$:
    \begin{enumerate}[(i)]
       \item If $x\notin A_\partial$, then
             $\operatorname{Rep}(i\otimes x)=\{(i,x)\}$.
        \item If $a_j\in   A_\partial$, then  $\operatorname{Rep}(i\otimes    a_j) \subseteq M \times A_\partial$. 
       \item If $a_j \in A_\partial$ and $A$ is non-degenerate, then $  \operatorname{Rep}(i\otimes a_j)=B_C$, where $C=[(i,j)]_{\overline R}$.
    \end{enumerate}
\end{lemma}

\begin{proof}
Every generating pair for $\approx_A$ consists of boundary points. Thus, no other point is related to a point of $M\times(A\setminus A_\partial)$. This proves (i) and (ii). When $ A$ is non-degenerate, (iii) follows from Lemma \ref{Lem:BoundaryClassCriterion}-(ii). The finiteness of $\operatorname{Rep}(u)$, for any $u \in \F A$ is now straightforward.
\end{proof}

\subsection{Non-degeneracy of a finite gluing specification}

To restrict the gluing functor to non-degenerate pointed objects, the selected outputs must remain distinct. The non-degeneracy of $\mathcal C$ as a $k$-pointed set is the precise requirement for this. Proposition~\ref{Prop_RNondegenerate_ModuleNondegenerate} in Appendix \ref{Subsec:RigorousModuleConstruction} identifies non-degeneracy of $\mathcal C$ with non-degeneracy of the associated module.

\begin{definition}[Non-degenerate finite gluing specification]
\label{Def_NondegeneracyR}
    Let $\Sigma=(M,D,R)$ be a finite $k$-pointed gluing specification.  We say that $\Sigma$ is \emph{non-degenerate} if $\mathcal C$ is non-degenerate:
    \begin{equation}
    \label{eqn:Non-degeneracy}
    \tag{ND}
        C_r \neq C_{r'} \quad \text{ for } \quad r,r' \in [k], \quad r\neq r'.
    \end{equation}
    Equivalently, $\eta:I\to\mathcal C$ is injective. In terms of the relation $R$, this means $(d_r,d_{r'})\notin \overline R$ whenever $r\ne r'$.
\end{definition}

\begin{remark}
    Non-degeneracy depends on $D$ as well as $R$: the same gluing relation can give distinct outputs for one choice of $D$ and coinciding outputs for another.
\end{remark}

The following proposition follows from Lemma~\ref{Lem:BoundaryClassCriterion} applied to the output labels, together with the functoriality established above.

\begin{proposition}[Non-degeneracy of $\F-$]
\label{Prop_EndofunctorOnSetND}
Let $\Sigma$ be a finite gluing specification. For every non-degenerate $k$-pointed set $A$,
\begin{equation}
\label{Lem_NonDegenerate_Set}
    i_r\otimes a_{j_r}=i_s\otimes a_{j_s} \quad\Longleftrightarrow\quad C_r=C_s \quad\text{for all }r,s\in[k].
\end{equation}
Consequently, the following are equivalent:
\begin{enumerate}[(i)]
    \item $\Sigma$ is non-degenerate;
    \item $\F A$ is non-degenerate for every non-degenerate $A$;
    \item $\F-$ restricts to an endofunctor on $\setnd$.
\end{enumerate}
\end{proposition}

For $C\in\mathcal C$, define its \emph{copy support} by
\begin{equation} \nonumber
   \operatorname{supp}(C) =\{i\in[\s]:(i,j)\in C\text{ for some }j\in[k]\}.
\end{equation}

\begin{definition}[Copy-separated outputs]
\label{Def_CopySeparatedOutputs}
    The finite gluing specification $\Sigma$ has \emph{copy-separated outputs} if
    \begin{equation}
        \operatorname{supp}(C_r)\cap\operatorname{supp}(C_{r'})=\emptyset \qquad \text{ for } r\neq r'.
    \label{eqn:Copy-separation}
    \tag{CS}
    \end{equation}
    
    Equivalently, no copy contains representatives of the $r$-th and $r'$-th outputs for two distinct indices $r\ne r'$.
\end{definition}

Equal output classes have the same nonempty copy support, contradicting copy-separation. Hence the following lemma holds. Also note that the converse fails because distinct classes can have representatives in the same copy. The metric role of copy-separation is established in Section~\ref{subsec_MetricAdmissibilityAndCopyseparation}.

\begin{lemma}
\label{Lem_CSImpliesNDButNotConversely}
Copy-separation implies non-degeneracy.
\end{lemma}

\begin{example}[Unit interval: separation]
\label{Ex_IntervalSeparation}
    The specification of Example~\ref{Ex_unitInterval} has boundary classes
    \begin{equation} \nonumber
        C_1=\{(1,1)\}, \qquad C_*=\{(1,2),(2,1)\}, \qquad C_2=\{(2,2)\}.
    \end{equation}
    The output classes $C_1$ and $C_2$ are distinct, with disjoint copy supports $\{1\}$ and $\{2\}$. Hence both non-degeneracy and copy-separation hold.
\end{example}

Our next goal is to extend $\F-$ on $\setnd$ to an endofunctor on $\mcnd$.

\subsection{Extension of $\M\otimes -$ to the metric setting}

Given a finite gluing specification $\Sigma=(M,D,R)$, $\w \in (0,1)$, and a $k$-pointed metric space $(A, d_A)$, we equip $M\times A$ with the metric given by
\begin{equation}
\nonumber
    \dm{A}\left((m,x),(n,y)\right)=\left\{
    \begin{array}{ll}
        \w \cdot d_A(x,y) & \text{ if }m=n \\
                            & \\
        1 & \hbox{ if } m\neq n
    \end{array}
    \right.
\end{equation}
Note that $\dm{A}$ on $M\times A$ is a one-bounded metric and hence $\left(\M\times A, \dm{A} \right) $ is a $k$-pointed metric space. We denote this $k$-pointed metric space by $\Mw A$. The metric on $\Mw A$ shrinks each copy of $A$ by a factor $\w \in (0, 1)$ and puts different copies at distance 1. In this section we allow any scaling factor, and the additional condition needed later for the canonical metric on the initial algebra is isolated later (cf. Definition \ref{Def_MetricAdmissibility}).

We shall equip $\F A$ with the quotient metric arising from the quotient map $\Mw A \to \F A$. We recall the construction. For $u,v\in \M \otimes A$, an \emph{admissible chain of $n$ steps} from $u$ to $v$ is a finite sequence
\begin{equation}
     \gamma : z_1,w_1,z_2,w_2,\ldots,z_n,w_n
\label{Eqn_paths}
\end{equation}
in $\M \times A$ such that
\begin{equation}
    z_1\in\operatorname{Rep}(u),\quad w_i\approx_A z_{i+1} \quad \text{ for } i\in[n-1],\quad \text{ and } \quad w_n\in\operatorname{Rep}(v).
\label{eqn:conditions on admisible chain}
\end{equation}

Its score is
\begin{equation*}
    \operatorname{sc}_A(\gamma) = \sum_{i=1}^n \dm{A}(z_i,w_i).
\label{Eqn_score}
\end{equation*}
We say such an admissible chain $\gamma$ is \emph{reduced} if $z_i\not\approx_A w_i$ for $1\leq i\leq n$ and $w_i\neq z_{i+1}$ for $1\leq i<n$. Define
\begin{equation}
    \dqm A\left(u,v\right) =\inf\{\operatorname{sc}_A(\gamma): \gamma\text{ is an admissible chain from }u\text{ to }v\}.
\label{eqn: quotient metric}
\end{equation}
The function $\dqm A$ is a pseudometric on $M\otimes A$. Because a one-term chain between arbitrary representatives of $u$ and $v$ has score at most one, $\dqm A(u,v)\leq1$.

\begin{definition}[Quotient metric on $\F -$]
\label{Def_MeticOnM_otimes_A}
    Given a finite gluing specification $\Sigma=(M,D,R)$, $\w \in (0,1)$ and a $k$-pointed metric space $(A, d_A)$, we let $\dqm{A}$ be the one-bounded quotient pseudometric on $\F A$ given by \eqref{eqn: quotient metric}. We write $\Fw A$ to denote the space $\left(\F A, \dqm{A}\right)$.
\end{definition}

We shall next show that for non-degenerate $A$, this pseudometric separates points and that every distance is attained by an admissible chain (cf. Proposition~\ref{Prop_FiniteBoundaryGraphFormula}). The metric construction holds for every $0<\w<1$.  However, in order to define the canonical metric on the initial algebra, we later require a metric admissibility condition.  The stronger condition $\w\geq \frac{1}{2}$ with copy-separated outputs is a convenient sufficient condition for this (cf. Proposition \ref{Prop_MetAdimisibilityAndCS}-(\ref{Prop_RhoHalfAdmissible})), but it is not necessary in all examples.

\subsubsection{The finite boundary graph}
\label{subsec:finite boundary graph}

Let $A$ be non-degenerate and fix $u,v\in\F A$.  Define a finite weighted graph $\Gamma_A(u,v)$ whose vertex set is the disjoint union
\begin{equation} \nonumber
   V=\{\mathsf s,\mathsf t\}\amalg\mathcal C.
\end{equation}
The auxiliary vertices $\mathsf s,\mathsf t$ represent the endpoints $u,v$, while the vertices in $\mathcal C$ represent the boundary gluing classes. Associate a finite subset $S_z\subseteq M \times A$ to each vertex by
\begin{equation} 
\nonumber
  S_{\mathsf s}=\operatorname{Rep}(u),   \qquad   S_{\mathsf t} =\operatorname{Rep}(v),   \qquad   S_C=B_C\quad(C\in\mathcal C).
\end{equation}
The sets $S_z\subseteq M\times A$ are finite by Lemma~\ref{Lem:FiniteQuotientFibres}. For distinct vertices $z,z'$, the edge weight is
\begin{equation}
    \ell(z,z')=  \min\Bigl\{   \bigl\{\rho d_A(x,y): (i,x)\in S_z,\ (i,y)\in S_{z'} \bigr\}\cup\{+\infty\}  \Bigr\}.
\label{eqn:edge Weight}
\end{equation}
An edge has finite length precisely when the two associated sets contain representatives lying in a common copy of $A$. The weight is the least cost of a metric step within a single scaled copy between representatives in the two associated sets. 

Let $L_A(u,v)$ denote the shortest-path length from $\mathsf s$ to $\mathsf t$, with $L_A(u,v)=+\infty$ if no path of finite weight exists. Since edge weights are non-negative, we can restrict to simple paths and there are only finitely many of those. Thus, whenever a finite-length path exists, the infimum over all finite-length walks is attained by a simple path.

\begin{proposition}[Finite graph formula and metric separation]
\label{Prop_FiniteBoundaryGraphFormula}
    Let $\Sigma$ be a gluing specification and $\w \in (0,1)$. Suppose $A$ is non-degenerate metric space. 
    \begin{enumerate}[(i)]
        \item For all $u,v\in \Fw A$,
            \begin{equation}
               \dqm{A}(u,v)=\min\{1,L_A(u,v)\}.
            \label{eqn: d = min(1, L_A) and Infimum Attained}
            \end{equation}
            Thus, the infimum in \eqref{eqn: quotient metric} defining $\dqm{A}(u,v)$ is attained.
        \item For $u,v\in \M\otimes A$, $\dqm{A}(u,v)$ is the minimum of the scores of the finite family of admissible  chains having $n<2+({\rho\Delta(A)})^{-1}$ metric steps, where every chain with $n>1$ is reduced and has score less than $1$.
        \item The quotient pseudometric $\dqm{A}$ is a one-bounded metric on $\F A$.   
        \item In particular, if $\Sigma$ is non-degenerate, then $\Fw A$ is also a non-degenerate $k$-pointed metric space.
    \end{enumerate}
\end{proposition}
\begin{proof}
For (i), we first prove $\min\{1,L_A(u,v)\}\leq\dqm{A}(u,v)$. If $\operatorname{sc}_A(\gamma)\geq1$, then $\min\{1,L_A(u,v)\}\leq\operatorname{sc}_A(\gamma)$ is immediate. Suppose therefore that its score is less than $1$.

Let $\gamma$ be any admissible chain from $u$ to $v$ (as in \eqref{Eqn_paths} and \eqref{eqn:conditions on admisible chain}) of  score less than $1$.  No metric step of $\gamma$ can join different copies, since every such step has cost $1$. We can simplify $\gamma$ without increasing its score.  If $w_i=z_{i+1}$ as points of $\M\times A$, then the two adjacent metric steps lie in the same copy and hence  can be replaced by a single metric step, by the triangle inequality in $A$. We also delete any metric step whose endpoints are $\approx_A$-equivalent, composing the adjacent equivalence jumps; the corresponding one-sided operation applies at either end of the chain. Repeating these operations terminates because the number of metric steps strictly decreases.

Unless $u=v$, the resulting chain has at least one metric step.  Every remaining equivalence jump is non-trivial.  Since all generators of $\approx_A$ are boundary pairs, the two points on either side of every such jump lie in sets $B_C$ for one and the same $C\in\mathcal C$.  The simplified chain therefore determines a graph walk from $\mathsf s$ to $\mathsf t$.  By \eqref{eqn:edge Weight}, each edge length of that walk is at most the length of the corresponding metric step.  Hence $ L_A(u,v)\leq\operatorname{sc}_A(\gamma)$. 

The same conclusion is trivial when $u=v$, because then $L_A(u,u)=0$.  Thus $\min\{1,L_A(u,v)\}$ is at most the score of every admissible chain.  Taking the infimum proves the inequality.

For the reverse inequality, first note that when $L_A(u,v)=+\infty$, the reverse follows immediately from $\dqm{A}(u,v)\leq1$. Otherwise, let $ \mathsf s=z_0,z_1,\ldots,z_r=\mathsf t $ be any finite-length graph path. For each edge $(z_{i-1},z_i)$, the minimum in \eqref{eqn:edge Weight} is attained because the sets involved are finite.  Choose
\begin{equation} 
\nonumber
   (p_i,x_i)\in S_{z_{i-1}} \; \text{ and } \; (p_i,y_i)\in S_{z_i} \quad \text{ with } \quad    \ell(z_{i-1},z_i)=\rho d_A(x_i,y_i).
\end{equation}
If $z_i=C\in\mathcal C$, then $(p_i,y_i)$ and $(p_{i+1},x_{i+1})$ both belong to $B_C$, and hence are $\approx_A$-equivalent, by Lemma~\ref{Lem:BoundaryClassCriterion}-(ii).  The chosen pairs consequently form an admissible chain from $u$ to $v$, of score exactly the length of the graph path.  It follows that
\begin{equation} 
\nonumber
   \dqm{A}(u,v)\leq L_A(u,v).
\end{equation}
Together with $\dqm{A}(u,v)\leq1$, this proves $\dqm{A}(u,v)\leq\min\{1,L_A(u,v)\}$. This completes the proof of \eqref{eqn: d = min(1, L_A) and Infimum Attained}.

If $\dqm{A}(u,v)=1$, every one-term admissible chain attains this value: its score is at most one and cannot be smaller than the infimum. If $\dqm{A}(u,v)<1$, then $L_A(u,v)=\dqm{A}(u,v)$, and a shortest graph path lifts to an admissible chain attaining this value. Thus the infimum in \eqref{eqn: quotient metric} is always attained.

For (ii), if a one-step chain attains $\dqm{A}(u,v)$, the assertion is immediate. Otherwise, $\dqm{A}(u,v)<1$ and (i) supplies a minimizing admissible chain. Reducing this chain preserves its minimum score and leaves $n\geq2$ metric steps, each within one copy. The $n-2$ internal metric steps join distinct boundary points, so
\begin{equation} 
\nonumber
    (n-2)\rho\Delta(A) \leq \rho\sum_{i=1}^{n}d_A(x_i,y_i) = \dqm{A}(u,v)<1.
\end{equation}
Hence $n<2+{(\rho\cdot \Delta(A)})^{-1}$.

Finally, the endpoint representatives belong to the finite sets $\operatorname{Rep}(u)$ and $\operatorname{Rep}(v)$, while all other representatives belong to the finite set $M\times A_\partial$. Since $n$ is bounded, the family of candidate chains is finite.

For (iii), suppose $\dqm{A}(u,v)=0$. By (ii), a chain from the specified finite family attains this value. Every reduced chain with more than one metric step has strictly positive score, so the attaining chain must be a one-step chain. Thus there exist $(p,x)\in\operatorname{Rep}(u)$ and $(q,y)\in\operatorname{Rep}(v)$ such that
\begin{equation} 
\nonumber
    0=\dqm{A}(u,v)=\dm{A}((p,x),(q,y)).
\end{equation}
Since $\dm{A}$ is a metric, $(p,x)=(q,y)$, and hence $u=v$. Therefore, $\dqm{A}$ is a one-bounded metric.

Proposition \ref{Prop_EndofunctorOnSetND} with (iii) implies (iv).
\end{proof}

\begin{example}[Unit interval: the quotient metric]
\label{Ex_IntervalQuotientMetric}
    Give $I=([2],1,2)$ the discrete metric. The boundary graph joins the two outputs through $C_*$:
\begin{equation} 
\nonumber
    C_1\xrightarrow{\rho}C_*\xrightarrow{\rho}C_2.
\end{equation}
This route costs $2\rho$, while a direct cross-copy step costs one. Consequently,
\begin{equation} 
\nonumber
    d_{\Fw I}(1\otimes1,2\otimes2)=\min\{1,2\rho\}.
\end{equation}
\end{example}

Metric separation does not require non-degeneracy of $\Sigma$, and its non-degeneracy guarantees the restriction to the non-degenerate metric category. We have the following estimate for $\Delta(\Fw A)$ when $A$ is non-degenerate and $1/2 \leq \w < 1$. 

\begin{lemma}[Separation estimate]
\label{Lem_DeltaOFMA}
    Assume that $\Sigma$ has copy-separated outputs and that $1/2\leq\w<1$. If $A$ is a non-degenerate $k$-pointed metric space, then
    \begin{equation} 
    \nonumber
       0<\Delta(A)\leq\Delta(\Fw A)\leq1.
    \end{equation}
    In particular, if $\Delta(A)=1$, then $\Delta(\Fw A)=1$.
\end{lemma}

\begin{proof}
The assertion is immediate for $k=1$, so suppose $k\geq2$.  Fix distinct outputs $r\neq r'$, and put
\begin{equation} 
\nonumber
   u_r=i_r\otimes a_{j_r}, \qquad u_{r'}=i_{r'}\otimes a_{j_{r'}}.
\end{equation}
By Proposition~\ref{Prop_FiniteBoundaryGraphFormula}, their distance is the smaller of $1$ and the length of a shortest simple path in the finite boundary graph from $C_r$ to $C_{r'}$.

A finite graph edge directly joining $C_r$ to $C_{r'}$ would be realized by representatives of these classes in the same copy, which violates \eqref{eqn:Copy-separation}.  Hence a same-copy boundary path from $C_r$ to $C_{r'}$ has at least two edges.  Every one of its edges joins two different boundary classes. Since $A$ is non-degenerate, the two corresponding distinguished points of $A$ are distinct, and the edge has length at least $\w\Delta(A)$.  Therefore every such path has length at least $2\w\Delta(A)\geq\Delta(A)$. The remaining candidate for the quotient distance is $1$, which is also at least $\Delta(A)$, because the metric of $A$ is one-bounded. Thus
\begin{equation} \nonumber
   d_{\Fw A}(u_r,u_{r'})\geq\Delta(A).
\end{equation}
Taking the minimum over $r\neq r'$ gives $\Delta(\Fw A)\geq\Delta(A)>0$.  The upper bound is the one-boundedness of the quotient metric.  If $\Delta(A)=1$, both bounds force $\Delta(\Fw A)=1$.
\end{proof}

\subsubsection{Local metric behavior and preservation of completeness}

For a non-degenerate $k$-pointed metric space $A$, give the finite boundary image $\mathcal C_A$ of \eqref{eqn:BoundaryImage} the metric inherited from $\Fw A$. It is canonically identified with $\mathcal C=\F I$ as a pointed set; its metric depends on $A$ and need not agree with the quotient metric on $\Fw I$. Since $\mathcal C_A$ is finite, it is closed in $\Fw A$. For $z\in\Fw A$, put
\begin{equation}
\nonumber
   r(z)=d_{\Fw A}(z,\mathcal C_A)=\min_{p\in \mathcal C_A}d_{\Fw A}(z,p).
\label{eqn:distance_to_gluing_locus}
\end{equation}

The following lemma records the local metric behavior away from the gluing locus $\mathcal C_A$, for $A \in \mcnd$.

\begin{lemma}
\label{lem:metric_behaviour_away_form_guling_locus}
    Let $A$ be a non-degenerate $k$-pointed metric space. If $z,z'\notin \mathcal C_A$, with unique representatives  $z=m\otimes x$ and $z'=m'\otimes x'$ in $M\otimes A$,  necessarily with $x, x'\notin A_{\partial}$,  then the following assertions hold.
    \begin{enumerate}[(i)]
        \item \label{eqn:bound_1} If $m\neq m'$ then   
        \begin{equation} 
        \nonumber
            d_{\Fw A}(z,z')\geq\min\{1,r(z)+r(z')\}.
        \end{equation}
        \item \label{eqn:bound_2} If $ m=m'$ and $d_{\Fw A}(z,z')<\min\{1,r(z)+r(z')\}$ then
        \begin{equation} 
        \nonumber
            d_{\Fw A}~(z,z')=\rho~d_A(x,x').
        \end{equation}
    \end{enumerate}
\end{lemma}

\begin{proof}
Choose a reduced minimizing chain for $d_{\Fw A}(z,z')$, guaranteed by Proposition~\ref{Prop_FiniteBoundaryGraphFormula}.  If its score is less than $1$, none of its metric steps can join two different copies, because every such step has cost $1$.  Moreover, every non-trivial equivalence jump in the chain occurs between boundary representatives and therefore passes through a point of $\mathcal C_A$.

Suppose that the chain contains at least one non-trivial equivalence jump. Let $p\in \mathcal C_A$ be the quotient point at its first such jump and let $p'\in \mathcal C_A$ be the quotient point at its last such jump.  Splitting the chain at these two jumps gives an initial subchain from $z$ to $p$ and a final subchain from $p'$ to $z'$.  If their scores are $s_0$ and $s_1$, respectively, then
\begin{equation} \nonumber
   s_0\geq d_{\Fw A}(z,p)\geq r(z), \qquad    s_1\geq d_{\Fw A}(p',z')\geq r(z').
\end{equation}
All intervening metric steps have non-negative length.  Hence, every such chain has score at least $ r(z)+r(z')$. Therefore, in this case
\begin{equation}
   d_{\Fw A}(z,z') \geq r(z)+r(z').
\label{eqn:bound_3}
\end{equation}

If $m\neq m'$, then a chain of score less than $1$ must contain a non-trivial equivalence jump: its metric steps remain within individual copies, so some equivalence jump is necessary to change the copy index. Estimate \eqref{eqn:bound_3}, together with the alternative that the minimum distance is $1$, proves \eqref{eqn:bound_1}.

Now suppose $m=m'$ and the strict inequality in \eqref{eqn:bound_2} holds.  A minimizing chain has score less than $1$, and \eqref{eqn:bound_3} shows that it cannot contain a non-trivial equivalence jump.  After the reductions in Proposition~\ref{Prop_FiniteBoundaryGraphFormula}, such a chain consists of a single metric step inside the $m$-copy.  Since $z$ and $z'$ have unique representatives, the score of this step is exactly $\rho d_A(x,x')$.  This proves \eqref{eqn:bound_2}, and hence the claim.
\end{proof}

A useful implication of the above result is the following lemma, which we use in Section~\ref{Sec:Completion} (see Appendix~\ref{App_FPreservesCompleteness} for a proof). 

\begin{lemma}
\label{Lem_FPreservesCompleteness}
    If $\Sigma$ is non-degenerate and $A$ is a complete non-degenerate $k$-pointed metric space, then so is $\Fw A$.
\end{lemma}

\subsubsection{Preservation Properties of $\Fw-$}
 To complete our goal of showing that $\Fw-$ extends to an endofunctor on $\mcnd$, we need to show that the induced maps preserve the following metric and topological properties. The following proposition is an extension of Lemma 2.3 of \cite{Sierpinski_gasket}.

\begin{proposition}
\label{Pro_PreservationProperties}
    Let $(\A, d_A)$ and $(\B, d_B)$ be non-degenerate $k$-pointed metric spaces and let $f:A \rightarrow B$ preserve distinguished points. For fixed $0<\rho<1$, the induced map $\Fw f$ has each of the following properties whenever $f$ does:
    \begin{enumerate}[(i)]
        \item \label{Lip} Lipschitz with Lipschitz constant $K\geq 1$
        \item \label{short} short map
        \item \label{isometry} isometric embedding
        \item \label{continuity} continuous
    \end{enumerate}
\end{proposition}

\begin{proof}
For \textup{(\ref{Lip})} (and \textup{(\ref{short})}, where $K=1$) apply $M\times f$ to a quotient chain. Gluing jumps remain gluing jumps. Each same-copy metric step increases by at most the factor $K$, while a between-copy step still has cost one, which is at most $K$ times its original cost. Taking infima over quotient chains gives
\begin{equation} 
\nonumber
    d_{\Fw B}\bigl((\Fw f)(u),(\Fw f)(v)\bigr) \leq K\,d_{\Fw A}(u,v).
\end{equation}

For (\ref{isometry}), suppose $f$ is an isometric embedding. Fix $u_1,u_2\in\Fw A$ and put $v_i=(\Fw f)(u_i)$. The graphs $\Gamma_A(u_1,u_2)$ and $\Gamma_B(v_1,v_2)$ have the same vertex labels. Their representative sets correspond under $M\times f$:
\begin{equation} 
\nonumber
    B_C^{\,B}=(M\times f)(B_C^{\,A}),\qquad \operatorname{Rep}_B(v_i)=(M\times f)(\operatorname{Rep}_A(u_i)).
\end{equation}
For the second equality, a non-trivial gluing zig-zag from an image point stays among distinguished points, all in $f(A)$; it pulls back along the injection by Proposition~\ref{Prop_ElementryPropertiesOfGluingConstruction}-(\ref{Lem_InjectionsReflectGluing}). Since $f$ preserves distances and copy indices, corresponding edges have equal weights, including $+\infty$. Thus $L_B(v_1,v_2)=L_A(u_1,u_2)$. Proposition~\ref{Prop_FiniteBoundaryGraphFormula} now gives 
\begin{equation} 
\nonumber
    d_{\Fw B}(v_1,v_2)=\min\{1,L_B(v_1,v_2)\} = d_{\Fw A}(u_1,u_2).
\end{equation}
Hence $\Fw f$ is an isometric embedding.

For (\ref{continuity}), suppose that $f:A\to B$ is continuous. Recall that $\mathcal C_A\subseteq\Fw A$ is finite and hence closed. We show that $\Fw f$ is continuous at $u\in\Fw A$.

First suppose $u\notin \mathcal C_A$. Then $u=m\otimes x$ has a unique representative, with $x\notin A_\partial$, and
\begin{equation} 
\nonumber
    r(u)=d_{\Fw A}(u,\mathcal C_A)>0.
\end{equation}
Given $\varepsilon>0$, continuity of $f$ at $x$ supplies $\delta'>0$ such that
\begin{equation} 
\nonumber
    d_A(x,y)<\delta' \quad\Longrightarrow\quad d_B(f(x),f(y))<\varepsilon/\rho.
\end{equation}
Choose $\delta=\frac12\min\{1,r(u),\rho\delta'\}$. If $d_{\Fw A}(u,v)<\delta$, then $v\notin \mathcal C_A$. By Lemma~\ref{lem:metric_behaviour_away_form_guling_locus}, $v$ cannot lie in a different copy: in that case its distance from $u$ would be at least $\min\{1,r(u)+r(v)\}\geq\min\{1,r(u)\}>\delta$. Thus $v=m\otimes y$ for some $y\notin A_\partial$. The same lemma now gives $d_{\Fw A}(u,v)=\rho d_A(x,y)$, because $d_{\Fw A}(u,v)<\min\{1,r(u)+r(v)\}$. Hence $d_A(x,y)<\delta'$, and a one-step chain in $\Fw B$ yields
\begin{equation} 
\nonumber
    d_{\Fw B}\bigl((\Fw f)(u),(\Fw f)(v)\bigr)  \leq \rho d_B(f(x),f(y)) <\varepsilon.
\end{equation}

Now suppose $u\in \mathcal C_A$. Since $A$ is non-degenerate, $\Delta(A)>0$. By continuity of $f$ at the finitely many distinguished points, there is $\delta'>0$ such that, for every $j\in[k]$,
\begin{equation} 
\nonumber
    d_A(a_j,y)<\delta'\quad \Longrightarrow \quad d_B(b_j,f(y))<\varepsilon/  \rho.
\end{equation}
Choose $\delta=\frac12\min\{1,\rho\Delta(A),\rho\delta'\}$. Let $d_{\Fw A}(u,v)<\delta$, and choose a minimizing chain as in Proposition~\ref{Prop_FiniteBoundaryGraphFormula}-(ii). This chain has only one metric step. Because the one-step chain has score less than $1$, its endpoints lie in the same copy. We may therefore write its representatives as $(m,a_j)\in\operatorname{Rep}(u)$ and $(m,y)\in\operatorname{Rep}(v)$, with $\rho d_A(a_j,y)=d_{\Fw A}(u,v)<\delta$. Thus $d_A(a_j,y)<\delta'$. Lifting this one-step chain along $f$ gives
\begin{equation} 
\nonumber
    d_{\Fw B}\bigl((\Fw f)(u),(\Fw f)(v)\bigr) \leq \rho d_B(b_j,f(y))  <\varepsilon.
\end{equation}
Hence $\Fw f$ is continuous at every point of $\Fw A$.
\end{proof}

The above observations lead to the main conclusion of this section, which follows from Proposition~\ref{Prop_EndofunctorOnSetND}, Proposition~\ref{Prop_FiniteBoundaryGraphFormula}  and  Proposition \ref{Pro_PreservationProperties}.

\begin{proposition}
\label{Prop:MetricEndofunctorNondegeneracy}
    Let $\Sigma=(M,D,R)$ be a finite gluing specification and $0<\rho<1$. Then, the metric gluing functor $\Fw-$ restricts to the category $\mcnd$ if and only if $\Sigma$ is non-degenerate.
\end{proposition}

\section{Initial algebra of $\F -$ and its canonical metric}
\label{Sec_InitialAlgebraOfF}
Fix a finite gluing specification $\Sigma=(M,D,R)$. Throughout this section and in the rest of the paper, unless otherwise mentioned, we assume that  $\Sigma$ is \textit{non-degenerate}. In particular, $\F -$ and $\Fw -$ restrict to endofunctors on $\setnd$ and $\mcnd$ respectively.

\subsection{Initial algebra of $\F -$ on $\set$ and $\setnd$}
\label{Subsec_InitialAlgebraOfFonSet}

The construction of the initial algebra follows from the standard results of \cite{AMM,Adamek74}. Recall that $I=([k],1,\ldots,k)$ is initial in both $\set$ and $\setnd$, and that $\mathcal C=\F I$ is the first stage after $I$. The map $\eta:I\to\mathcal C$ is given by $\eta(r)=C_r=i_r\otimes j_r$. Write $\Fi{0}-$ for the identity functor and $\Fi{n}-$ for the $n$-fold iterate of $\F-$. The initial chain arising from the initial object is
\begin{equation}
\label{inichain}
     I\xrightarrow{\eta}\F I \xrightarrow{\F\eta}\Fi{2}I  \xrightarrow{\Fi{2}\eta}\cdots.
\end{equation}
Since $\Sigma$ is non-degenerate, all its stages are non-degenerate and $\eta$ is injective. Proposition~\ref{Prop_ElementryPropertiesOfGluingConstruction}-(\ref{Lem_oneone}) then makes every connecting map injective.

For $r\geq s\geq0$, let $\eta_{r,s}:\Fi{s}I\to\Fi{r}I$ be the canonical map along this chain, with $\eta_{s,s}=\operatorname{id}_{\Fi{s}I}$. The colimit of \eqref{inichain} in $\set$ is 
\begin{equation} 
\nonumber  
    G=\left(\coprod_{n\geq0}\Fi{n}I\right)\big/\cong,
\end{equation}
where, for $x\in\Fi{r}I$ and $y\in\Fi{s}I$,
\begin{equation} 
\nonumber
    x\cong y \quad \Longleftrightarrow \quad \eta_{t,r}(x)=\eta_{t,s}(y) \quad \text{for some } t\geq r,s.
\end{equation}
Write $[x]$ for the class of $x$. The canonical maps $\pi_n:\Fi{n}I\to G$, $\pi_n(x)=[x]$, are injective and satisfy $\pi_r\circ\eta_{r,s}=\pi_s$. Thus $G$ is non-degenerate, with distinguished points $[1],\ldots,[k]$. Every element has a finite-tensor representative 
\begin{equation} 
\nonumber
    [m_1\otimes\cdots\otimes m_n\otimes j],  \qquad m_i\in M,\quad j\in[k],
\end{equation}
with $n=0$ giving $[j]$.

Given the above construction, the following proposition follows from Ad\'amek's initial-algebra theorem \cite{Adamek74}. We note that $\F-$ preserves the colimit of \eqref{inichain}. Since $G$ is non-degenerate and $\setnd$ is a full subcategory of $\set$, initiality holds in both categories. The structure map is an isomorphism by Lambek's lemma \cite{Lam}.

\begin{proposition}[Initial algebra]
\label{Prop_initialAlgebraOnSet}
The map
\begin{equation} 
\nonumber
    g:\F G\longrightarrow G, \qquad g(m\otimes[x])=[m\otimes x],
\end{equation}
makes $(G,g)$ the initial $\F-$algebra in both $\set$ and $\setnd$, and $g$ is an isomorphism. For an $\F-$algebra $(X,e)$ in $\set$, put 
\begin{equation} 
\nonumber
    e_0=\operatorname{id}_X, \qquad e_{n+1}=e\circ(\F e_n).
\end{equation}
The unique algebra morphism $f:G\to X$ is given by 
\begin{equation} 
\nonumber
    f([m_1\otimes\cdots\otimes m_n\otimes j]) = e_n(m_1\otimes\cdots\otimes m_n\otimes x_j),
\end{equation}
where $x_j$ is the $j$-th distinguished point of $X$. Equivalently, $f$ is the unique pointed map satisfying
\begin{equation}
\label{Diag_initialSquare}
    f\circ g=e\circ(\F f).
\end{equation}
\end{proposition}

\subsection{Metric Admissibility and Copy-Separation}
\label{subsec_MetricAdmissibilityAndCopyseparation}

Equip $I=([k],1,\ldots,k)$ with the discrete metric $d_I$, so that $d_I(r,s)=1$ for $r\ne s$ and $d_I(r,r)=0$. Thus $I$ is non-degenerate and $\Delta(I)=1$, including the convention for $k=1$. We also write $I$ for this pointed metric space. The canonical metric on the initial algebra is obtained by making all maps in the initial chain isometric. Since $\Fw-$ preserves isometric embeddings, it is enough to check its first connecting map. 

\begin{definition}[Metric admissibility]
\label{Def_MetricAdmissibility}
    Assume $0<\w<1$ and that the fixed gluing specification $\Sigma$ is non-degenerate. We say that the scaled specification $(\Sigma,\w)$ is \emph{metrically admissible} if the canonical map
    \begin{equation}
        \eta:I\longrightarrow \Fw I \quad \text{  is an isometric embedding, }
    \label{eqn:MetricAdmisibility}
    \tag{MA}
    \end{equation}
     where $I$ is equipped with the discrete metric.
\end{definition}

\begin{proposition}[Metric admissibility and copy-separation]
    \label{Prop_MetAdimisibilityAndCS}
        Let $\Sigma$ be a non-degenerate finite gluing specification and $0<\rho<1$.
    \begin{enumerate}[(i)]
        \item \label{Lem_AdmissibilityImpliesCopySeparation}
            Metric admissibility of $(\Sigma,\rho)$ implies copy-separation.
        \item \label{Prop_RhoHalfAdmissible}
            If $\frac12\leq\rho<1$, then $(\Sigma,\rho)$ is metrically admissible if and only if $\Sigma$ has copy-separated outputs.
    \end{enumerate}
\end{proposition}

\begin{proof}
Suppose distinct output classes have representatives in one copy, say $(i,j)$ and $(i,j')$. Then
\begin{equation} 
\nonumber
    d_{\Fw I}\bigl(\eta(r),\eta(s)\bigr) \leq\rho d_I(j,j')\leq\rho<1=d_I(r,s).
\end{equation}
Thus $\eta$ is not isometric, proving (i). Conversely, copy-separation implies non-degeneracy, and for $\rho\geq\frac12$ the separation estimate (Lemma~\ref{Lem_DeltaOFMA}) gives $\Delta(\Fw I)=1$. Hence $\eta:I\to\Fw I$ preserves all distances, proving (ii).
\end{proof}

Thus, for $0<\rho<1$,
\begin{equation} 
\nonumber
    \mathrm{(MA)} \Longrightarrow \mathrm{(CS)} \Longrightarrow \mathrm{(ND)},
\end{equation}
with $\mathrm{(MA)}\Longleftrightarrow\mathrm{(CS)}$ when $\frac12\leq\rho<1$. The categorical scale $\rho$ need not equal an IFS contraction ratio; Section~\ref{Sec:IFS} gives the additional comparison hypotheses.

Write $\Fwi{i} -$ for the $i$-fold composition of $\Fw -$. Assume that $(\Sigma,\w)$ is metrically admissible.  By Proposition \ref{Pro_PreservationProperties}, all the maps in the metric initial chain 
\begin{equation}
    I~\xrightarrow{\eta}~ \Fw I ~\xrightarrow{\Fw \eta}~ \Fwi{2} I ~\xrightarrow{\Fwi{2} \eta} ~ \Fwi{3} I \rightarrow \cdots \rightarrow \Fwi{i} I ~\xrightarrow{\Fwi{i} \eta}\cdots~~~
\label{inichain-metric} 
\tag{\ref{inichain}'}
\end{equation}
are also isometric embeddings. The carrier set $G$ of the initial algebra in the set setting can be given the quotient metric. As all the maps in the initial chain \eqref{inichain-metric} are isometric embeddings, the quotient metric is a metric and can be computed as follows. 

For $x\in\Fwi{r}I$ and $y\in\Fwi{s}I$, choose $t\geq r,s$ and set
\begin{equation} 
\nonumber
    d_G([x],[y]) = d_{\Fwi{t}I}\bigl(\eta_{t,r}(x),\eta_{t,s}(y)\bigr).
\end{equation}
Since the connecting maps are isometric embeddings, this is independent of the chosen stage and representatives.

Thus $(G, [1],\cdots [k])$, with the metric $d_G$, is a non-degenerate $k$-pointed metric space with $\Delta(G)=1$. By abuse of notation, we shall also write  $G$ to mean the metric space $(G\,,\,d_G)$, and it should be understood from the context that $G$ refers to the underlying set  when the metric structure is not relevant. Moreover, $\Fw G$ is also a non-degenerate $k$-pointed metric space. To lift $G$ to an $\Fw-$ algebra over $\mcnd$, we need to show that $g:\Fw G \rightarrow G$ is continuous. In fact, $g$ is an isometric isomorphism. Consequently, $(G,g)$ is an algebra in the metric category.

\begin{proposition}
\label{Prop_GasAMetricAlgebra}
    Under metric admissibility of $\Sigma$, $g:\Fw G\to G$ is an isometric isomorphism. In particular, $(G,g)$ is an $\Fw-$algebra in $\mcnd$.
\end{proposition}

\begin{proof}
The map $g$ is bijective by initiality. Write $\pi_n:\Fwi{n}I\to G$ for the canonical isometric embeddings. Every pair $u,v\in\Fw G$ can be written at a common stage as $u=(\Fw\pi_n)(u')$ and $v=(\Fw\pi_n)(v')$. Since $\Fw\pi_n$ is isometric and $g\circ\Fw\pi_n=\pi_{n+1}$,
\begin{equation} 
\nonumber
    d_G(g(u),g(v)) = d_{\Fwi{n+1}I}(u',v') = d_{\Fw G}(u,v).
\end{equation}
Thus $g$ is an isometric isomorphism.
\end{proof}

\begin{example}[Unit interval: the initial algebra]
\label{Ex_IntervalInitialAlgebra}
    Since $\eta(1)=1\otimes1$ and $\eta(2)=2\otimes2$, Example~\ref{Ex_IntervalQuotientMetric} shows that metric admissibility holds exactly when $\frac12\leq\rho<1$. Take $\rho=\frac12$. 
    
    We represent copy $i\in[2]$ by the binary digit $i-1$. The $n$-th stage identifies with $  D_n=\left\{\frac{j}{2^n}:j=0,\ldots,2^n\right\}$, and the connecting maps become inclusions. As each new stage joins two half-size grids end-to-end, its quotient metric is Euclidean. Thus the initial algebra $G$ that consists of the finite binary terms modulo $\overline R$ can be identified with the dyadic rationals
        \begin{equation} \nonumber
        \mathbb D = \left\{ \frac{m}{2^n} \;\middle|\; n\geq0,\ 0\leq m\leq2^n \right\} \subseteq [0,1].
    \end{equation}
    Under this identification, the canonical metric on $G$ agrees with the Euclidean metric on $\mathbb D$. 

    Transporting the algebra structure of $G$ along the identification $G\cong\mathbb D$ gives the initial algebra $g:\F\mathbb D\longrightarrow\mathbb D$ of $\F-$ on $\setndk{2}$, given by
    \begin{equation} \nonumber
       g(1\otimes q)=\frac q2, \qquad g(2\otimes q)=\frac q2+\frac12.
    \end{equation}
\end{example}

\section{Completion of the Initial Algebra}
\label{Sec:Completion}

We assume $(\Sigma,\w)$ is metrically admissible. In particular $G$ is equipped with its canonical metric $d_G$ and $(G, g)$ can  be considered as a $\Fw-$algebra over $\mcnd$. We apply Cauchy completion functorially to this metric algebra. Thus, by ``completion of the initial algebra" we mean the completion of this set-theoretic initial algebra with respect to its canonical metric. 

Even though $(G, g)$ is the initial algebra for $\F -$ on $\setnd$, it is important to note that $(G,g)$ need not be initial in the metric setting, i.e., that of $\Fw -$. For example, in the case of Sierpiński  gasket it is $G_d$, the same set with the discrete metric, that forms the initial algebra in the metric setting (see Lemma 4.1 of \cite{Sierpinski_gasket}). But, completion with respect to the discrete metric is itself and does not produce the final coalgebra. Thus, to  emphasize further, it is the completion with respect to the canonical metric that is used below. 

To this end, we first recall the Cauchy completion as a  functorial construction on the Lipschitz subcategory of $\mcnd$. Completion commutes with $\Fw -$ up to a natural isometric isomorphism. This gives the appropriate notion of completing an algebra whose structure map is Lipschitz.  We then apply this construction to the algebra $(G,g)$.

\subsection{Cauchy completion}
For a non-degenerate $k$-pointed metric space $\X$, let $C(X)$ be its completion. Denote the class of a Cauchy sequence $(a_i)$  by $[a_i]$ and the class of the constant sequence  of $a\in X$ by $[\bar{a}]$. We make $C(X)$ a $k$-pointed metric space by taking the $r$-th distinguished point of $C(X)$ to be $[\bar{x_r}]$. Note that $C(X)$ is also non-degenerate. $X$ is isometrically embedded in $C(X)$ via the canonical pointed map $x \mapsto [\bar{x}]$.  
For a Lipschitz pointed map $f:X\to Y$,
\begin{equation} \nonumber
    C(f):C(X)\longrightarrow C(Y),\qquad C(f)([a_i])=[f(a_i)],
\end{equation}
is the extension of $f$ to $C(X)$, and hence it is also pointed and has the same Lipschitz constant. Completion is therefore functorial on the Lipschitz subcategory of $k$-pointed metric spaces and it preserves isometric embeddings. The isometric embedding of $X$ into $C(X)$ gives rise to an isometric embedding of $\Fw X$ into $\Fw C(X)$ (by Proposition \ref{Pro_PreservationProperties}).  

Given $X$ in $\mcnd$, one can first Cauchy complete $X$ and apply $\Fw -$ to obtain $\Fw C(X)$. On the other hand, one can first consider $\Fw X$ and take the completion of it, $C(\Fw X)$. These two constructions are naturally isometrically isomorphic on the Lipschitz subcategory. The following result generalizes \cite[Lemma~13]{Bhat} and is recalled  in \cite[Proposition~3.5]{Sierpinski_gasket}.

\begin{proposition}[Lemma~13 of \cite{Bhat}]
\label{prop_BetaIsNatural}
    For $X$ in $\mcnd$, the map 
    \begin{equation}
    \label{Def_beta}
    \begin{aligned}
        \beta_X:\Fw C(X)  &  \longrightarrow C(\Fw X),\\
        m\otimes[a_n]     &  \longmapsto[m\otimes a_n]
    \end{aligned}
    \end{equation}
    is an isometric isomorphism. The collection of isometries $\{\beta_X\}$ forms a natural isomorphism from $(\Fw-)\circ C$ to $ C\circ(\Fw-)$ on the Lipschitz subcategory of $\mcnd$.
\end{proposition}

\begin{proof}
Let $\iota_X:X\to C(X)$ be the canonical isometric embedding. By Proposition~\ref{Pro_PreservationProperties}, $\Fw\iota_X:\Fw X\to\Fw C(X)$ is an isometric embedding. Its image is dense: if $a_n\to a$ in $C(X)$, with $a_n\in X$, then $m\otimes a_n\to m\otimes a$. By Lemma~\ref{Lem_FPreservesCompleteness}, $\Fw C(X)$ is complete. Consequently, $\Fw\iota_X$ extends uniquely to an isometric isomorphism $\gamma_X:C(\Fw X)\to\Fw C(X)$. Its inverse is the displayed map $\beta_X$. Naturality follows by checking the relevant identity on the dense copy of $\Fw X$ and then extending by continuity.
\end{proof}

Write 
\begin{equation}
\nonumber
    \gamma_X=\beta_X^{-1}:C(\Fw X)\to\Fw C(X).
\end{equation}
Note that the above result is sufficient for the completion of the initial algebra $(G,g)$, as the structure map $g:\Fw G\rightarrow G$ is an isometric isomorphism.  The final coalgebra theorem itself is still stated in the category $\mcnd$ of continuous maps.

\subsection{Completion of the initial algebra of $\F-$}
\label{Sec_CompletionOfTheInitialAlgebra}

\begin{definition}[Completion of an algebra]
\label{Def_CompletionOfAlgebra}
For an $\Fw-$algebra $(X,e)$ in $\mcnd$ with Lipschitz structure map, its completion is $(C(X),\tilde e)$, where
\begin{equation} 
\nonumber
    \tilde e=C(e)\circ\beta_X:\Fw C(X)\longrightarrow C(X).
\end{equation}
\end{definition}

Given an algebra $(X,e)$ as in the above definition, the canonical embedding $\iota_X: X \to C(X)$ is an algebra morphism:
\begin{equation}
\label{Eqn_embedingAlgebraToCompletion}
\begin{tikzcd}
    \Fw X \arrow[r,"e"] \arrow[d,hookrightarrow,"\Fw\iota_X"'] & X \arrow[d,hookrightarrow,"\iota_X"]\\
    \Fw C(X) \arrow[r,"\tilde e"'] & C(X).
\end{tikzcd}
\end{equation}

Applying this construction to the isometric isomorphism $g:\Fw G\to G$ gives
\begin{equation} 
\nonumber
    \G=C(G),\qquad \tilde g=C(g)\circ\beta_G:\Fw\G\longrightarrow\G.
\end{equation}
Since $C(g)$ and $\beta_G$ are isometric isomorphisms, so is $\tilde g$. Its inverse 
\begin{equation} 
\nonumber
    \g=\gamma_G\circ C(g^{-1}):\G\longrightarrow\Fw\G
\end{equation}
therefore makes $(\G,\g)$ an $\Fw-$coalgebra in $\mcnd$. 

\begin{example}[Unit interval: completion]
\label{Ex_IntervalCompletion}
    Since $\mathbb D$ is dense in $[0,1]$, the identification of $G$ with $\mathbb D$ extends uniquely to an isometric isomorphism $ eval:\G=C(G)\to [0,1]$, given by the binary evaluation map, where $a_n\in\{0,1\}$ represents copy $a_n+1$, 
    \begin{equation} 
    \nonumber
         eval\big([a_0a_1a_2\cdots]\big) = \sum_{n=0}^{\infty}\frac{a_n}{2^{n+1}}.
    \end{equation}
    Under this identification, the completed algebra structure is given by 
    \begin{equation} \nonumber
        \tilde g(i\otimes x)=(i-1+x)/2, \quad i\in[2].
    \end{equation}
\end{example}

\section{Final Coalgebra of $\F-$ and $\Fw -$}
\label{Sec:FinalCoalgebra}
We continue to assume $(\Sigma,\w)$ is metrically admissible. In this section, we show that $(\G, \g)$ is the final coalgebra of $\F-$ on $\setnd$. The canonical metric on $\G$ will play a role in determining the mediating morphism, even in the set setting. Note that we do \underline{not} assert the finality of $(\G,\g)$ for $\F-$ on all of $\set$.  We also prove that the mediating morphism associated with every continuous $\Fw-$coalgebra is continuous.  

Throughout we shall view $\G$ as a complete $k$-pointed metric space containing $G$ as a dense subset. This inclusion is a map of $k$-pointed metric spaces, which induces a map of $\Fw -$ algebras, as in diagram (\ref{Eqn_embedingAlgebraToCompletion}). In particular, the distinguished points of $\G$ are the same as those of $G$. Moreover, the metric on $\G$ is the extension of that on $G$ (i.e., $d_G$) to $\G$ via limits.

\subsection{Final Coalgebra of $\F-$ on $\setnd$}

We can now complete the proof of Theorem \ref{Thm_FinalCoalgebraInSet}. We adapt the contraction-mapping argument of \cite[proofs of Theorems~4 and~14]{Bhat}, and then recover the finite-stage description of the mediating morphism. Recall the following result first.

For a word $\mathbf m=(m_1,\ldots,m_n)\in M^n$, write $\mathbf m\otimes y=m_1\otimes\cdots\otimes m_n\otimes y$.  Given $X$ in $\mcnd$,  $n\geq1$, $\mathbf m\in M^n$, and $y,z\in X$, the bound 
\begin{equation} 
\label{Eqn_Mi_inequality}
    d_{\Fwi{n}X} \bigl(\mathbf m\otimes y,\mathbf m\otimes z\bigr) \leq \w^n d_X(y,z) \leq \w^n.
\end{equation}
generalizes Lemma 15 of \cite{Bhat} and is proved similarly. 

\begin{proof}[Proof of Theorem~\ref{Thm_FinalCoalgebraInSet}]
The initial-algebra assertion follows from Proposition~\ref{Prop_initialAlgebraOnSet}. It remains to prove finality.

Equip the space $\G^X$ of all functions with the supremum metric
\begin{equation} 
\nonumber
    d_\infty(u,v)=\sup_{x\in X}d_\G(u(x),v(x)).
\end{equation}
Because $\G$ is complete and one-bounded, this metric makes $\G^X$ complete. Since $X$ is non-degenerate, 
\begin{equation} 
\nonumber
    \mathcal H=\operatorname{Hom}_{\setnd}(X,\G)
\end{equation}
is a nonempty closed subset of $\G^X$, as it consists of the functions satisfying $u(x_j)=[j]$ for every $j\in[k]$. In particular, $\mathcal H$ is complete.

Let $e:X\to\F X$ be a coalgebra in $\setnd$. Define
\begin{equation} 
\nonumber
    T:\mathcal H\longrightarrow\mathcal H, \qquad T(u)=\tilde g\circ(\F u)\circ e.
\end{equation}
The composition preserves distinguished points. For $u,v\in\mathcal H$, choose a representative $e(x)=m\otimes y$. Since $\tilde g$ is an isometry,
\begin{equation} 
\nonumber
\begin{aligned}
    d_\G\bigl(T(u)(x),T(v)(x)\bigr)  & = d_{\Fw\G} \bigl(m\otimes u(y),m\otimes v(y)\bigr)\\
        &   \leq \w\,d_\G(u(y),v(y))\\
        &   \leq \w\,d_\infty(u,v).
\end{aligned}
\end{equation}
Taking suprema shows that $T$ is a contraction with constant $\w<1$. By the Banach fixed-point theorem, $T$ has a unique fixed point $f\in\mathcal H$. The fixed-point equation 
\begin{equation} 
\nonumber
    f=\tilde g\circ(\F f)\circ e
\end{equation}
is equivalent to the commutativity of
\begin{equation}
\label{Eqn_comSquareMead}
\begin{tikzcd}
    X \arrow[r,"e"] \arrow[d,"f"'] & \F X \arrow[d,"\F f"]\\
    \G \arrow[r,"\g"']
        & \F\G .
\end{tikzcd}
\end{equation}
Thus $f$ is the unique coalgebra morphism from $(X,e)$ to $(\G,\g)$, proving finality in $\setnd$.
\end{proof}

We now recover the finite-stage description of this map. For an $\F-$coalgebra $(X,e)$, define $e_n:X\to\Fi{n}X$ recursively by
\begin{equation} 
\nonumber
    e_0=\operatorname{id}_X, \qquad e_{n+1}=(\Fi{n}e)\circ e_n.
\end{equation}

\begin{lemma}
\label{Lem_MediatingMapApproximation}
    Let $(X,e)$ be an $\F-$coalgebra in $\setnd$, and let $f:X\to\G$ be its unique coalgebra morphism. Fix $r\in[k]$. For each $n\geq1$ and $x\in X$, choose any representative
    \begin{equation} 
    \nonumber
        e_n(x)=m_1\otimes\cdots\otimes m_n\otimes y = \mathbf m\otimes y,
    \end{equation}
    and put
    \begin{equation} 
    \nonumber
        \theta_n(x)=[\mathbf m\otimes r]\in G.
    \end{equation}
    Then
    \begin{equation} 
    \nonumber
        d_\G\bigl(f(x),\theta_n(x)\bigr)\leq\w^n.
    \end{equation}
    Consequently, $\theta_n\to f$ uniformly in $\G$, independent of the representatives and the choice of $r$. 
    
    In particular, independent of all auxiliary choices, for any $x\in X$,
    \begin{equation} 
    \nonumber
        f(x)=\lim_{i\rightarrow \infty}\theta_i(x).
    \end{equation}
\end{lemma}

\begin{proof}
Put $\tilde g_0=\operatorname{id}_\G$ and
\begin{equation} 
\nonumber
    \tilde g_{n+1}=\tilde g\circ(\Fw\tilde g_n).
\end{equation}
Each $\tilde g_n:\Fwi{n}\G\to\G$ is an isometric isomorphism.

Iterating the coalgebra identity for $f$, and using that $G\hookrightarrow\G$ is an algebra morphism, gives
\begin{equation} 
\nonumber
    f(x)=\tilde g_n(\mathbf m\otimes f(y)), \qquad \theta_n(x)=\tilde g_n(\mathbf m\otimes[r]).
\end{equation}
Hence \eqref{Eqn_Mi_inequality} yields
\begin{equation} 
\nonumber
    d_\G\bigl(f(x),\theta_n(x)\bigr) = d_{\Fwi{n}\G} \bigl(\mathbf m\otimes f(y), \mathbf m\otimes[r]\bigr) \leq\w^n.
\end{equation}
The bound is uniform and applies to every choice of representatives and distinguished point.
\end{proof}

\subsection{Final Coalgebra of $\Fw-$ on $\mcnd$}
\label{Sec:FinalCoalgebraMetric}
We now show that $(\G\,,\,\g)$ is also the final coalgebra of $\Fw-$ on $\mcnd$. Let $(X\,,\,e)$ be an $\Fw-$coalgebra. Thus, $X$ is a non-degenerate $k$-pointed metric space and $e:X\rightarrow \Fw X$ is continuous. By forgetting the metric structure, we obtain from Theorem \ref{Thm_FinalCoalgebraInSet} a unique map $f:X\rightarrow \G$ such that the diagram in \eqref{Eqn_comSquareMead} commutes.  To prove finality in the metric category, it remains to show that this mediating morphism is continuous.

\begin{lemma}[Finite-stage estimate for the mediating map]
\label{Lem_MediatingMapEstimate}
    Let $(X,e)$ be an $\Fw-$coalgebra on $\mcnd$, and let $f:X\rightarrow\G$ be the mediating map obtained in Theorem~\ref{Thm_FinalCoalgebraInSet}.  Then, for all $x,y\in X$ and every $i\geq1$,
    \begin{equation}
    \label{Eqn_MediatingMapEstimate}
        d_\G(f(x),f(y)) \leq 2\w^i + \frac{1}{\Delta(X)} \left(2\w^i + d_{\Fwi{i} X}\bigl(e_i(x),e_i(y)\bigr) \right).
    \end{equation}
\end{lemma}

\begin{proof}
Given $X$, consider $X_\partial$ with the induced metric and the induced $k$-pointed structure.  The inclusion $\iota:X_\partial\hookrightarrow X$ is an isometric embedding, and hence so is
\begin{equation} 
\nonumber
    \Fwi{i}\iota:\Fwi{i}X_\partial\longrightarrow \Fwi{i}X
\end{equation}
for every $i\geq 0$.  Moreover, with the convention that $\Delta=1$ in the one-pointed case, we have $\Delta(X_\partial)=\Delta(X)$.  Since $X$ is non-degenerate and all metrics are one-bounded, $0<\Delta(X)\leq 1$, and hence $1/\Delta(X)\geq 1$.

Define $\alpha:X_\partial\rightarrow I$ by $\alpha(x_j)=j$.   If $j\neq j'$, then
\begin{equation} 
\nonumber
    d_I(\alpha(x_j),\alpha(x_{j'}))=1\leq \frac{1}{\Delta(X)}d_X(x_j,x_{j'}).
\end{equation}
Thus $\alpha$ is Lipschitz with constant $1/\Delta(X)$.  By Proposition \ref{Pro_PreservationProperties}, each
\begin{equation} 
\nonumber
    \Fwi{i}\alpha:\Fwi{i}X_\partial\longrightarrow \Fwi{i}I
\end{equation}
is also Lipschitz with constant $1/\Delta(X)$.

Fix $x,y\in X$ and $i\geq1$.  Choose representatives $ e_i(x)=\mathbf{n}\otimes a$ and $e_i(y)=\mathbf{n}'\otimes b$, where $\mathbf{n}=n_1\otimes\cdots\otimes n_i$, $\mathbf{n}'=n'_1\otimes\cdots\otimes n'_i$, and $a,b\in X$.  
Then,
\begin{equation} 
\nonumber
    \theta_i(x)=[\mathbf{n}\otimes r], \quad \text{ and } \quad     \theta_i(y)=[\mathbf{n}'\otimes r].
\end{equation}

By Lemma~\ref{Lem_MediatingMapApproximation},
\begin{equation}
\label{eqn:estimate_FSE_bound_1}
    d_\G(f(z),\theta_i(z))\leq\w^i \qquad \text{for } z\in\{x,y\}.
\end{equation}

Using the isometric embedding $\Fwi{i}I\rightarrow G$, the Lipschitz estimate for $\Fwi{i}\alpha$, and the isometric embedding $\Fwi{i}X_\partial\rightarrow \Fwi{i}X$, we obtain
\begin{equation} 
\nonumber
\begin{aligned}
    d_G(\theta_i(x),\theta_i(y))
        & = d_{\Fwi{i}I}(\mathbf{n}\otimes r,\mathbf{n}'\otimes r) \\
        & \leq \frac{1}{\Delta(X)} d_{\Fwi{i}X_\partial}(\mathbf{n}\otimes x_r,\mathbf{n}'\otimes x_r) \\
        & = \frac{1}{\Delta(X)}       d_{\Fwi{i}X}(\mathbf{n}\otimes x_r,\mathbf{n}'\otimes x_r) \\
        & \leq \frac{1}{\Delta(X)}\Bigl( d_{\Fwi{i}X}(\mathbf{n}\otimes x_r,\mathbf{n}\otimes a)       +d_{\Fwi{i}X}(\mathbf{n}\otimes a,\mathbf{n}'\otimes b) \\
        & \hspace{96pt} +d_{\Fwi{i}X}(\mathbf{n}'\otimes b,\mathbf{n}'\otimes x_r) \Bigr) \\
        & \leq \frac{1}{\Delta(X)} \left(2\w^i+d_{\Fwi{i}X}(e_i(x),e_i(y))\right).
\end{aligned}
\end{equation}
Therefore, combining this inequality with \eqref{eqn:estimate_FSE_bound_1} and the triangle inequality gives
\begin{equation} 
\nonumber
\begin{aligned}
    d_\G(f(x),f(y)) & \leq d_\G(f(x),\theta_i(x)) +d_G(\theta_i(x),\theta_i(y))        +d_\G(\theta_i(y),f(y)) \\
        & \leq 2\w^i+ \frac{1}{\Delta(X)} \left( 2\w^i+d_{\Fwi{i}X}(e_i(x),e_i(y)) \right),
\end{aligned}
\end{equation}
which is \eqref{Eqn_MediatingMapEstimate}.
\end{proof}

We now complete the proof of Theorem \ref{Thm_FinalCoalgebraInMet}. 

\begin{proof}[Proof of Theorem \ref{Thm_FinalCoalgebraInMet} ]
Fix $x\in X$ and $\epsilon>0$.  Choose $q\geq1$ so large that
\begin{equation} 
\nonumber
    2\w^q\left(1+\frac{1}{\Delta(X)}\right)<\frac{\epsilon}{2}.
\end{equation}
The iterate $e_q:X\rightarrow\Fwi{q}X$ is continuous, because it is a composition of $e$ and its images under the endofunctor $\Fw-$.  Hence there is $\delta>0$ such that
\begin{equation} 
\nonumber
    d_X(x,y)<\delta  \quad \Longrightarrow \quad  d_{\Fwi{q}X}\bigl(e_q(x),e_q(y)\bigr)  <\frac{\epsilon\Delta(X)}{2}.
\end{equation}
Lemma~\ref{Lem_MediatingMapEstimate} now gives
\begin{equation} 
\nonumber
    d_\G(f(x),f(y))<\frac{\epsilon}{2}+\frac{\epsilon}{2} = \epsilon.
\end{equation}
Thus $f$ is continuous. Its uniqueness as a continuous coalgebra morphism is the set-theoretic uniqueness established in Theorem~\ref{Thm_FinalCoalgebraInSet}.  Therefore $(\G,\g)$ is final in $\mcnd$.
\end{proof}

\begin{example}[Unit interval: the final coalgebra]
\label{Ex_IntervalFinalCoalgebra}
    Under the isometry of Example~\ref{Ex_IntervalCompletion}, the inverse of the completed algebra is 
    \begin{equation} \nonumber
        \gamma:[0,1]\longrightarrow \Fww{\frac12}[0,1], \qquad \gamma(x)=
        \begin{cases}
            1\otimes2x, & 0\leq x\leq\frac12,\\[1mm]
            2\otimes(2x-1), & \frac12\leq x\leq1.
        \end{cases}
    \end{equation}
    At $x=\frac12$ the expressions agree because $1\otimes1=2\otimes0$. By Theorems~\ref{Thm_FinalCoalgebraInSet} and~\ref{Thm_FinalCoalgebraInMet}, this coalgebra is final in $\mcndk{2}$ and, after forgetting the metric, in $\setndk{2}$. Note that this realizes the self-similarity of the unit interval as described in the introduction.
\end{example}

 It follows from Proposition \ref{Pro_PreservationProperties} that $\Fw-$ restricts to an endofunctor on the subcategory of $\mcnd$ consisting of short maps. Since $\g$ is an isometry, $(\G,\g)$ is also an $\Fw-$coalgebra in this subcategory. However, the mediating map from an arbitrary non-degenerate short-map coalgebra need not be short unless the distinguished points of the domain are normalized. The following proposition gives the precise estimate.

\begin{proposition}
\label{Prop_MediatingMapLipschitz}
    Assume that $(\Sigma,\w)$ is metrically admissible. 
    \begin{enumerate}[(i)]
        \item Let $(X,e)$ be an $\Fw-$coalgebra in $\mcnd$ with a short     structure map $e:X\to \Fw X$. Then the mediating map $f:X\rightarrow \G$ obtained in Theorem~\ref{Thm_FinalCoalgebraInSet} is Lipschitz with constant $1/\Delta(X)$. In particular, if $\Delta(X)=1$, then $f$ is short.
        \item The functor $\Fw-$ restricts to the subcategory of $\mcnd$ consisting of those non-degenerate $k$-pointed metric spaces $X$ such that $\Delta(X)=1$ and short maps. Moreover, $(\G,\g)$ is the final coalgebra of $\Fw-$ on this subcategory.
    \end{enumerate} 
\end{proposition}

\begin{proof}
For (i), note that every iteration of $e$ is short, as $e$ is short and $\Fw-$ preserves short maps. In particular, their compositions 
\begin{equation} 
\nonumber
    e_q:X\longrightarrow \Fwi{q}X
\end{equation}
are also short. Hence
\begin{equation} 
\nonumber
    d_{\Fwi{q}X}(e_q(x),e_q(y))\leq d_X(x,y).
\end{equation}
Lemma~\ref{Lem_MediatingMapEstimate} therefore gives, for every $q\geq1$,
\begin{equation} 
\nonumber
    d_\G(f(x),f(y)) \leq  2\w^q\left(1+\frac{1}{\Delta(X)}\right)  + \frac{1}{\Delta(X)}d_X(x,y).
\end{equation}
Letting $q\to\infty$ yields
\begin{equation} 
\nonumber
    d_\G(f(x),f(y))\leq \frac{1}{\Delta(X)}d_X(x,y).
\end{equation}
Therefore $f$ is Lipschitz with constant $1/\Delta(X)$. If $\Delta(X)=1$, then $f$ is short.

For (ii), assume $\Delta(X)=1$. Then the canonical pointed map $I\rightarrow X$ is an isometric embedding.  Proposition~\ref{Pro_PreservationProperties} therefore makes $\Fw I\rightarrow\Fw X$ an isometric embedding.  Metric admissibility gives $\Delta(\Fw I)=1$, and hence $\Delta(\Fw X)=1$. This shows the restriction of $\Fw-$ to the subcategory. By metric admissibility, the canonical inclusion $I\rightarrow \G$ is also an isometry and pointed. Hence $\Delta(\G)=1$. Moreover, $\g:\G\rightarrow \Fw\G$ is an isometry. Therefore $(\G,\g)$ is an $\Fw-$coalgebra in this subcategory.

Let $(X,e)$ be any $\Fw-$coalgebra in this subcategory. The unique underlying set-theoretic coalgebra morphism $f:X\longrightarrow \G$ is short because $\Delta(X)=1$ and by Proposition~\ref{Prop_MediatingMapLipschitz}. Therefore $f$ is the required unique coalgebra morphism. Hence $(\G,\g)$ is final.
\end{proof}

Note that uniform rescaling of the metric to force $\Delta=1$ may in tern change the one-boundedness.

\section{Application to Iterated Function Systems}
\label{Sec:IFS}
Let $(X,d_X)$ be a complete one-bounded metric space and $(\f_1,\ldots,\f_k)$, where $k>1$, be an IFS with common Lipschitz constant $\lambda<1$. Let $\alpha_i$ be the fixed point of $\f_i$ and $O_0=\{\alpha_1,\ldots,\alpha_k\}$. We make $O_0$ a $k$-pointed metric space by choosing the $r$-th distinguished point to be the $r$-th fixed point $\alpha_r$. Put $\Delta = \Delta(O_0)$. Let $O$ be the orbit of fixed points under finite compositions of the branches. By Hutchinson's theorem~\cite{Hut}, the IFS has a unique nonempty compact attractor $\mathcal F$ such that
\begin{equation} 
\nonumber
    \mathcal F=\bigcup_{i=1}^k\f_i(\mathcal F)=\overline O.
\end{equation}
We construct the associated gluing specification and compare $O\subseteq\mathcal F$ with $G\subseteq C(G)$, proving Theorem~\ref{Thm_F_isTheFinalCoalgebra}.

\subsection{The associated IFS gluing specification}
\label{subsec_AssociatedIFSGluingSpec}

\begin{definition}[Gluing specification of an IFS]
\label{Def_RelationForIFS}
An IFS $(\f_1,\ldots,\f_k)$ with fixed points $(\alpha_1,\cdots, \alpha_k)$ determines a finite $k$-pointed gluing specification $\Sigma=(M,D,R)$ by
\begin{equation} 
\nonumber
    M=[k],\qquad D=((1,1),\ldots,(k,k)),
\end{equation}
and
\begin{equation} 
\nonumber
    (i,j)R(i',j') \quad \Longleftrightarrow \quad \f_i(\alpha_j)=\f_{i'}(\alpha_{j'}).
\end{equation}
In particular, $R$ is the kernel relation of $(i,j)\mapsto\f_i(\alpha_j)$, and is already an equivalence relation.
\end{definition}

The associated gluing specification is the fixed-point coincidence relation of \cite[Example~10.11]{Lein}, where the attractor is identified as a universal solution using \cite[Corollary~9.7]{Lein}. Here we compare the orbit algebra with the initial algebra and then complete it; metric admissibility and a bi-Lipschitz comparison then identify $C(G)$ with $\mathcal F$. 

\subsubsection{Non-degeneracy}

Throughout, we assume that the fixed points are pairwise distinct. This assumption has two purposes. First, under this assumption, $O_0=\{\alpha_1,\ldots,\alpha_k\}$ and its orbit $O$ become non-degenerate $k$-pointed spaces (cf. Section \ref{subsec_OrbithAndAlgebraStructure}). Also, this assumption makes the gluing specification $\Sigma$ non-degenerate.  

\begin{lemma}
\label{lem:ND-distinct-fixed-points}
    The associated gluing specification is non-degenerate if and only if the fixed points are pairwise distinct. 

    Under the standing assumption that the fixed points are distinct, the set and metric functors therefore restrict to the non-degenerate categories.
\end{lemma}

\begin{proof}
Since $R=\overline R$ and $\f_r(\alpha_r)=\alpha_r$, the output classes $C_r=[(r,r)]_R$ satisfy
\begin{equation} 
\nonumber
    C_r=C_{r'} \iff \f_r(\alpha_r)=\f_{r'}(\alpha_{r'}) \iff \alpha_r=\alpha_{r'}.
\end{equation}
Now apply Propositions~\ref{Prop_EndofunctorOnSetND} and~\ref{Prop:MetricEndofunctorNondegeneracy}.
\end{proof}

\subsubsection{Output-isolation and copy-separation}

\begin{definition}[Output-isolated IFS]
\label{Def:IFSOutputIsolation}
    The IFS is \emph{output-isolated} if
    \begin{equation}
        \f_i(O_0)\cap O_0=\{\alpha_i\}\qquad \text{ for all } \quad i\in[k].
     \label{eqn:Output-isolation}
     \tag{OI}
    \end{equation}
\end{definition}

\begin{lemma}
\label{lem:IFS-OI-iff-CS}
Under the distinct-fixed-points assumption, the IFS is output-isolated if and only if its specification $\Sigma$ has copy-separated outputs.

In particular, if $\lambda<\Delta$, then the IFS is output-isolated and hence its specification has copy-separated outputs. In addition, if $\w \geq 1/2$ then $\Sigma$ is also metrically admissible.
\end{lemma}

\begin{proof}
The class $C_r$ contains $(r,r)$, and $i\in\operatorname{supp}(C_r)$ exactly when $\alpha_r\in\f_i(O_0)$. Output-isolation therefore says $\operatorname{supp}(C_r)=\{r\}$ for every $r$. Since each support already contains $r$, this is equivalent to pairwise disjointness of the output supports.

Now, the assumption $\f_i(\alpha_j)=\alpha_r$ with $r\ne i$ yields \begin{equation} 
\nonumber
    \Delta\leq d_X(\alpha_i,\alpha_r) = d_X(\f_i(\alpha_i),\f_i(\alpha_j))  \leq\lambda d_X(\alpha_i,\alpha_j)\leq\lambda.
\end{equation}
Thus, $\lambda < \Delta$ implies output isolation, and hence the copy separation. Metric admissibility for $\w \geq 1/2$ follows from Proposition \ref{Prop_MetAdimisibilityAndCS}.
\end{proof}

Consequently, for $\frac12\leq\rho<1$, output-isolation, copy-separation and metric admissibility are equivalent by Proposition~\ref{Prop_MetAdimisibilityAndCS}; the inequality $\lambda<\Delta$ is a sufficient check. For $0<\rho<\frac12$, copy-separation remains necessary and admissibility must be checked directly. These conditions do not imply orbit-algebra injectivity that we require later.

\subsection{The orbit and its algebra structure}
\label{subsec_OrbithAndAlgebraStructure}
Starting from the fixed points $O_0$, define
\begin{equation} 
\nonumber
    O_r=\bigcup_{i=1}^k\f_i(O_{r-1})\quad(r\geq1) \quad \text{ and } \qquad O=\bigcup_{r\geq0}O_r.
\end{equation}
Because $\alpha_i=\f_i(\alpha_i)$, the finite stages are ascending, and both $O_r$ and $O$ inherit the metric and distinguished points $(\alpha_1,\ldots,\alpha_k)$. They are non-degenerate (as fixed points are distinct) and $\overline O=\mathcal F$.

\begin{proposition}[Orbit algebra and finite stages]
\label{prop_OrbithAlgebraAndFiniteStages}
    Put $K=\max\{1,\lambda/\rho\}$. The formulas
    \begin{equation} 
    \nonumber
        {}_r\f:\Fw O_{r-1}\to O_r, \qquad \f:\Fw O\to O,\qquad i\otimes x\longmapsto\f_i(x)
    \end{equation}
    define surjective $K$-Lipschitz pointed maps for every $r\geq1$. In particular, $(O,\f)$ is an algebra in $\mcnd$, and these maps are short when $\lambda\leq\rho$. Moreover, $\f$ is injective if and only if every $ {}_r\f$ is injective.
\end{proposition}

\begin{proof}
The defining relation $R$ identifies boundary representatives with the same image under the branch maps. Hence the displayed maps are well-defined and pointed. They are surjective because $O_r=\bigcup_i\f_i(O_{r-1})$ and $O=\bigcup_i\f_i(O)$.

Within a copy, a scaled metric step has length $\rho d_X(x,y)$ and its image has length at most $\lambda d_X(x,y)$. A metric step between copies has length one, while its image has distance at most one. Gluing jumps have equal images. Summing along quotient chains proves the bound $K$, uniformly in $r$ and for $O$.

The finite-stage maps are restrictions of $\f$ along the injective maps $\F O_{r-1}\to\F O$. Thus injectivity of $\f$ implies stagewise injectivity. Conversely, any two representatives in $\F O$ belong to a common finite stage, where injectivity of $ {}_r\f$ applies.
\end{proof}

To identify the completed orbit with a final coalgebra, the algebra map $\f:\Fw O\rightarrow O$ needs to be inverted, so that its inverse supplies a coalgebra structure. Surjectivity leaves injectivity as the extra condition needed for this. As per the above proposition, injectivity of $\f$ amounts to the injectivity of $ {}_r\f $ for all $r\in \mathbb N$. The injectivity of orbit algebra must be checked separately from non-degeneracy, copy-separation, and metric admissibility.

\subsection{Set-theoretic comparison of the orbit algebra with the initial algebra}

Initiality of $(G, g)$ gives a unique algebra morphism $f:G\to O$, satisfying $f\circ g=\f\circ(\F f)$. Let $u:I\to O_0$ be the canonical pointed bijection and $\pi_r:\Fi{r}I\to G$ the colimit maps.

\begin{proposition}[Comparison with the orbit algebra]
\label{pro_IsoO_o_and_G_inSet}
    The unique algebra morphism $f:G\to O$ is given by 
    \begin{equation} 
    \nonumber
        f([i_1\otimes\cdots\otimes i_r\otimes p]) = (\f_{i_1} \circ \cdots \circ \f_{i_r})(\alpha_p).
    \end{equation}
    If $\f:\F O\to O$ is injective, then $f:G\rightarrow O$     is an isomorphism of $\F-$ algebras $(G,g)$ and $(O,\f)$ (on $\setnd$). In particular, $(O,\f)$ is initial in $\setnd$.
\end{proposition}

\begin{proof}
The formula follows from the universal property of $(G,g)$ and the evaluation formula in Proposition~\ref{Prop_initialAlgebraOnSet}.

Assume that $\f$ is injective. Each $ {}_r\f$ is then bijective. Define
\begin{align}
    {}_0\tf & = \operatorname{id}_{O_0},\qquad {}_1\tf={}_1\f,
    \label{Def_bijection_MnO_o_and_O_r}\\
    {}_{r+1}\tf&={}_{r+1}\f\circ(\F\,{}_r\tf) \qquad(r\geq1).
    \label{Def_bijection_MnO_o_and_O_r_GeneralCase}
\end{align}
Inductively, $ {}_r\tf:\Fi{r}O_0\to O_r$ is bijective, and
\begin{equation} 
\nonumber
    f\circ\pi_r={}_r\tf\circ(\Fi{r}u).
\end{equation}
Hence $f$ restricts to a bijection from each stage $\pi_r(\Fi{r}I)$ onto $O_r$. Any two elements of $G$ belong to a common stage, so $f$ is injective; its image contains every $O_r$, so it is surjective. A bijective algebra morphism is an algebra isomorphism.
\end{proof}

The inverse $h=f^{-1}$ restricts on $O_r$ to
\begin{equation} 
\nonumber
    {}_rh=\pi_r\circ(\Fi{r}u^{-1})\circ{}_r\tf^{-1}.
\end{equation}
These restrictions are compatible because they are restrictions of the same inverse $h$.

\subsection{Metric comparison with the initial algebra}

Assume throughout that $(\Sigma,\rho)$ is metrically admissible. Then $g:\Fw G\to G$ is isometric, and Proposition~\ref{pro_IsoO_o_and_G_inSet} gives an algebra morphism $f:G\to O$.

\begin{lemma}
\label{Lem_f_short}
    If $\lambda\leq\rho$, the comparison map $f:G\to O$ is short.
\end{lemma}

\begin{proof}
The pointed bijection $u:I\to O_0$ is short. Each ${}_r\f$ is short by Proposition~\ref{prop_OrbithAlgebraAndFiniteStages}, so the recursion \eqref{Def_bijection_MnO_o_and_O_r_GeneralCase} makes ${}_r\tf\circ(\Fwi{r}u)$ short. These are the restrictions of $f$ to the isometric initial stages. Any two points of $G$ belong to a common stage, proving the claim.
\end{proof}

\begin{proposition}[Topological comparison]
\label{pro_Iso_G_O_metricSetting}
    Suppose $(\Sigma,\rho)$ is metrically admissible and the orbit algebra $\f:\Fw O\to O$ is injective.
    \begin{enumerate}[(i)]
        \item If $\f^{-1}$ is continuous, then $h=f^{-1}$ is continuous.
        \item If $f$ is continuous, then $h$ is continuous if and only if $\f^{-1}$ is continuous. Under these equivalent conditions, $f$ is an algebra isomorphism in $\mcnd$.
    \end{enumerate}
\end{proposition}

\begin{proof}
If $\f^{-1}$ is continuous, $(O,\f^{-1})$ is a continuous coalgebra. With $\iota:G\hookrightarrow\G$, the map $\iota\circ h$ is its set-theoretic coalgebra morphism into $(\G,\g)$, and is continuous by Theorem~\ref{Thm_FinalCoalgebraInMet}. Since $\iota$ is an embedding, $h$ is continuous. Conversely, when both $f$ and $h$ are continuous, the algebra identity gives
\begin{equation} 
\nonumber
    \f^{-1}=(\Fw f)\circ g^{-1}\circ h,
\end{equation}
which is continuous by Proposition~\ref{Pro_PreservationProperties}.
\end{proof}

\begin{remark}[One-step inverse bounds]
\label{prop_Lipofphi_inverse}
The inverse $\f^{-1}$ is $\kappa$-Lipschitz if and only if every one-step inverse ${}_r\f^{-1}:O_r\to\Fw O_{r-1}$ has that bound. Indeed, the inclusions of orbit stages, and their images under $\Fw-$, are isometric, and every pair of orbit points lies in a common stage. Such a bound makes $\f$ bi-Lipschitz. It does not by itself provide a uniform bound on the \emph{iterated} coding inverses ${}_r\tf^{-1}$, since bounds can accumulate under composition. A uniform estimate controls $h=f^{-1}$, as follows.
\end{remark}

\begin{proposition}[Iterated inverse bounds]
\label{prop_Lipof_tphi_inverse}
    Suppose $(\Sigma,\rho)$ is metrically admissible, the orbit algebra $\f:\Fw O\to O$ is injective and ${}_r\tf^{-1}:O_r\to\Fwi{r}O_0$ is $\kappa$-Lipschitz for every $r\geq1$, with the same constant $\kappa>0$. Then:
    \begin{enumerate}[(i)]
        \item $h=f^{-1}:O\to G$ is $\kappa\Delta^{-1}$-Lipschitz.
        \item If $f$ is continuous, then $\f$ is a homeomorphism and $f$ is
            an algebra isomorphism in $\mcnd$.
        \item If $\lambda\leq\rho$, then both $\f$ and $f$ are bi-Lipschitz.
    \end{enumerate}
\end{proposition}

\begin{proof}
The map $u^{-1}:O_0\to I$ is $\Delta^{-1}$-Lipschitz, with $\Delta^{-1}\geq1$. Proposition~\ref{Pro_PreservationProperties} preserves this bound under every iterate. Thus 
\begin{equation} 
\nonumber
    h|_{O_r}=\pi_r\circ(\Fwi{r}u^{-1})\circ{}_r\tf^{-1}
\end{equation}
is $\kappa\Delta^{-1}$-Lipschitz. Passing to a common finite stage proves (i). If $f$ is continuous, (i) makes it a homeomorphism, and the identities
\begin{equation} 
\nonumber
    \f=f\circ g\circ(\Fw h), \qquad \f^{-1}=(\Fw f)\circ g^{-1}\circ h
\end{equation}
give (ii). When $\lambda\leq\rho$, Lemma~\ref{Lem_f_short} makes $f$ short, and the same identities make both $\f$ and its inverse Lipschitz, proving (iii).
\end{proof}

\subsection{Completion and identification with the final coalgebra}
Suppose, in addition to metric admissibility, that the orbit algebra $\f:\Fw O\to O$ is injective. Since $X$ is complete and $O$ is dense in the attractor $\mathcal F$, the inclusion $O\hookrightarrow X$ extends to the canonical isometry 
\begin{equation} 
\nonumber
    \kappa_O:C(O)\longrightarrow\overline O=\mathcal F, \qquad \kappa_O([x_n])=\lim_{n\to\infty}x_n.
\end{equation}
The orbit algebra is Lipschitz by Proposition~\ref{prop_OrbithAlgebraAndFiniteStages}. Its completion is therefore
\begin{equation} 
\nonumber
    \widehat{\f} = C(\f)\circ\beta_O: M\otimes_\rho C(O)\longrightarrow C(O),
\end{equation}
and transport along $\kappa_O$ is
\begin{equation} 
\nonumber
    \widetilde{\f} = \kappa_O\circ \widehat{\f}\circ (\Fw\kappa_O^{-1}):\Fw\mathcal F \longrightarrow\mathcal F.
\end{equation}
In particular, $\widetilde{\f}(i\otimes x)=\f_i(x)$ on the attractor, and hence completion extends the original IFS branches. With these definitions, we can now prove Theorem \ref{Thm_F_isTheFinalCoalgebra}.

\begin{proof}[Proof of Theorem~\ref{Thm_F_isTheFinalCoalgebra}.]
Proposition~\ref{pro_IsoO_o_and_G_inSet} gives the algebra isomorphism $f:G\to O$; its assumed bi-Lipschitz property proves (i). It extends to a bi-Lipschitz isomorphism $C(f):C(G)\to C(O)$. Put $\widetilde f=\kappa_O\circ C(f)$. Naturality of $\beta$ and $f\circ g=\f\circ(\Fw f)$ give the single transport identity
\begin{equation} 
\nonumber
    \ft\circ\widetilde g = \widetilde{\f}\circ(\Fw\widetilde f).
\end{equation}
Since $\widetilde f$ and $\widetilde g$ are isomorphisms, $\widetilde{\f}$ is an isomorphism and $\widetilde f$ identifies their inverse coalgebras. Finality now transfers from $(\G,\g)$ by Theorems~\ref{Thm_FinalCoalgebraInMet} and~\ref{Thm_FinalCoalgebraInSet}, proving (ii) and (iii).
\end{proof}

\begin{example}[Unit interval: the IFS comparison]
\label{Ex_unitInterval_IFS}
    On $[0,1]$, take $\f_i(x)=(i-1+x)/2$ for $i\in[2]$. The fixed points are $\alpha_1=0$ and $\alpha_2=1$, so $\Delta=1$, $\lambda=1/2$, and the associated specification is that of Example~\ref{Ex_unitInterval}. Lemma \ref{lem:ND-distinct-fixed-points} recovers the non-degeneracy of $\Sigma$. Its output supports are disjoint; hence $\rho=1/2$ is admissible.

    Because $\lambda<\Delta$, Lemma~\ref{lem:IFS-OI-iff-CS} gives output isolation and therefore copy-separation, which we claimed earlier. Metric admissibility follows from Proposition \ref{Prop_MetAdimisibilityAndCS}-(\ref{Prop_RhoHalfAdmissible}). By Lemma~\ref{Lem_f_short} the comparison $f:G\to O$ is short. 

    The orbit $O$ consists of the dyadic points, with closure $[0,1]$. The branch images meet only at $1/2$, exactly the prescribed boundary identification. The orbit algebra is therefore injective. It is also isometric: within a half-interval distances are halved; across the two halves a minimizing chain passes through $1/2$, and its score is the Euclidean distance. Induction identifies every initial stage isometrically with its dyadic image, so $f:G\to O$ is an isometry. Proposition~\ref{pro_IsoO_o_and_G_inSet} identifies $(O,\f)$ with the initial algebra $(G,g)$ in the set-theoretic setting and Proposition~\ref{pro_Iso_G_O_metricSetting} identifies $G$ and $O$ homeomorphically. Theorem~\ref{Thm_F_isTheFinalCoalgebra} now recovers the final coalgebra in Example~\ref{Ex_IntervalFinalCoalgebra} for $\Fww{\frac12}-$ on $\mcndk{2}$, and for $\F-$ on $\setndk{2}$.   

    In Example~\ref{Ex_Nested_Radicals_IFS}, we give an application of the universality of $[0,1]$ to nested radicals of 2.
\end{example}

\paragraph{IFS verification checklist.}  The conditions have distinct roles: non-degeneracy keeps the outputs distinct, admissibility supplies the canonical metric on $G$, orbit-algebra injectivity identifies $G$ with $O$, and quantitative comparison identifies their completions. The parameters $\lambda$, $\Delta$, and $\rho$ control, respectively, contraction in the ambient space, separation of the fixed points, and the scale of the categorical metric. Table~\ref{tab:IFS_conditions} collects the comparisons.

\begin{table}[htbp]
\centering
\small
\renewcommand{\arraystretch}{1.2}
\begin{tabularx}
 {\linewidth}{@{}
 >{\raggedright\arraybackslash}p{0.18\linewidth}
 >{\raggedright\arraybackslash}X
 >{\raggedright\arraybackslash}p{0.25\linewidth}@{}}
\toprule
Requirement & Criterion or sufficient check & Reference\\
\midrule
\midrule
Non-degeneracy
    & The fixed points are pairwise distinct.
    & Lemma~\ref{lem:ND-distinct-fixed-points}\\
\midrule
Metric admissibility
    & Verify that $\eta:I\to\Fw I$ is isometric. For $\rho\geq\frac12$, this is equivalent to copy-separation; $\lambda<\Delta$ suffices. For $\rho<\frac12$, check directly.
    & Proposition~\ref{Prop_MetAdimisibilityAndCS} \par Lemma \ref{lem:IFS-OI-iff-CS}\\
\midrule
Orbit-algebra injectivity
    & Every finite-stage map ${}_r\f$ is injective.
    & Proposition~\ref{prop_OrbithAlgebraAndFiniteStages}\\
\midrule
Bi-Lipschitz comparison
    & Direct estimates on $f$ and $f^{-1}$, or $\lambda\leq\rho$ and a uniform Lipschitz bound on all ${}_r\tf^{-1}$.
    & Proposition~\ref{prop_Lipof_tphi_inverse}\\
\bottomrule
\end{tabularx}
\caption{IFS verification checklist.}
\label{tab:IFS_conditions}
\end{table}

\section{Examples}
\label{Sec:Examples}
The unit interval has served as a running example. For the standard spaces below, we summarize the established coalgebraic presentations and explain their IFS realizations. The final two examples use interval finality to evaluate nested radicals and exhibit the obstruction caused by segment overlaps.

\begin{example}[Cantor set $\mathcal C$ and the binary streams]
\label{Ex_CantorSet}
    We recover the Cantor set by completing the algebra of eventually constant binary streams and compare the resulting prefix metric with its Euclidean metric. This recovers the metric version of the Cantor set and its head--tail presentation of binary streams (see \cite[Example~10.1]{Lein}).

    \subparagraph{Presentation}  Let $k=2$, the copy set $M=[2]$, $D=((1,1),(2,2))$, and $R$ the equality relation. Then,
    \begin{equation} 
    \nonumber
        \F(A,a_1,a_2) = \bigl([2]\times A,(1,a_1),(2,a_2)\bigr).
    \end{equation}
    There is no gluing between the two copies. The initial algebra is the set $G$ of eventually constant binary streams, pointed by $0^\omega$ and $1^\omega$; its algebra map is $g(i,s)=(i-1)s$, where copy $i\in[2]$ prepends the binary digit $i-1$. For every $0<\rho<1$, the specification is metrically admissible and the canonical metric is
    \begin{equation} 
    \nonumber
        d_\rho(a,b)=\rho^{\ell(a,b)}\quad(a\ne b), \qquad d_\rho(a,a)=0,
    \end{equation}
    where $\ell(a,b)$ is the length of the longest common prefix.     Its completion is all of $2^\omega$, with final coalgebra $\g(a_0a_1\cdots)=(a_0+1)\otimes(a_1a_2\cdots)$.
    
    \subparagraph{IFS realization} The Cantor set is the attractor of the IFS on $[0,1]$ given by 
    \begin{equation} 
    \nonumber
        \f_1(x)=x/3,\qquad \f_2(x)=(x+2)/3,
    \end{equation}
    whose fixed points are $0$ and $1$. From this, we recover the above non-degenerate gluing specification. They have disjoint branch images. The orbit $O$ of the two fixed points is the set of endpoints of the finite Cantor approximations. The ternary evaluation
    \begin{equation}
    \label{eqn:TirnaryEvaluation_Cantor}
        f(a_0a_1\cdots)=\sum_{n=0}^\infty\frac{2a_n}{3^{n+1}}
    \end{equation}
    is the unique algebra morphism that identifies $G$ with the orbit $O$. Since the branches are injective and have disjoint images, the orbit algebra is injective.     For distinct streams with common prefix of length $\ell$,
    \begin{equation} 
    \nonumber
        \frac13\,d_{1/3}(a,b) \leq |f(a)-f(b)| \leq d_{1/3}(a,b).
    \end{equation}
    Thus $\rho=1/3$ gives a bi-Lipschitz comparison, and     Theorem~\ref{Thm_F_isTheFinalCoalgebra} identifies the attractor     with the completed final coalgebra $(2^\omega, \g)$. The transported coalgebra structure on the attractor $\mathcal C$ is given by $\g_{\mathcal C}(\f_i(x))=i\otimes x$.
    
    All the prefix metrics $d_\rho$, $0<\rho<1$, induce the     product topology on $2^\omega$. Thus ternary evaluation     \eqref{eqn:TirnaryEvaluation_Cantor} extends $f$ to a homeomorphism
    \begin{equation} 
    \nonumber
        \widehat f:(2^\omega,d_\rho)\longrightarrow\mathcal C
    \end{equation}
    for every $0<\rho<1$. The canonical comparison is bi-Lipschitz precisely when $\rho=1/3$: the streams $0^\omega$ and $0^n1^\omega$ have prefix distance $\rho^n$, while their images have Euclidean distance $3^{-n}$.
\end{example}

\begin{example}[Sierpiński gasket $\Ss$]
\label{Ex_Sierpinski}
    For the standard three-copy gluing functor, the set-theoretic initial algebra is the union of the finite gasket vertex sets \cite[Proposition 4.2]{Sierpinski_gasket}. Its canonical completion is bi-Lipschitz equivalent to the Euclidean gasket \cite{Bhat,Sierpinski_gasket}. The continuous final coalgebra presentation is established in \cite[Theorem 4.4]{Sierpinski_gasket}.
    
    \subparagraph{IFS realization} On the closed equilateral triangle with vertices $\alpha_1=(1/2,\sqrt3/2)$,  $\alpha_2=(0,0)$ and $\alpha_3=(1,0)$,  consider the IFS given by 
    \begin{equation} 
    \nonumber
        \f_i(z)=\frac{(z+\alpha_i)}{2}, \quad \text{ where } \; i\in[3].
    \end{equation}
    The fixed points are precisely $\alpha_1$, $\alpha_2$ and $\alpha_3$. The associated specification recovers that of \cite{Sierpinski_gasket}. It has $k=3$, copy set $M=[3]$, $D=((1,1),(2,2),(3,3))$, and $R$ generated by
    \begin{equation} 
    \nonumber
        (1,2)R(2,1),\qquad (1,3)R(3,1),\qquad (2,3)R(3,2).
    \end{equation}
    This specification is non-degenerate and output classes are copy-separated. This also makes $\rho=1/2$ metrically admissible (by Proposition~\ref{Prop_MetAdimisibilityAndCS}). The first-level pieces satisfy
    \begin{equation} 
    \nonumber
        \f_i(\Ss)\cap\f_j(\Ss) = \{\f_i(\alpha_j)\} = \{\f_j(\alpha_i)\}      \qquad (i\ne j).
    \end{equation}
    Since each branch is injective and $R$ identifies precisely these pairs of boundary representatives, the orbit algebra is injective. The orbit $O$ is the union of the finite gasket vertex sets and is dense in $\Ss$.
    
    For $\rho=1/2$, by \cite[Theorem~19]{Bhat} the comparison $f:G\to O$ satisfies, 
    \begin{equation} 
    \nonumber
        \frac{\sqrt3}{4}\,d_G(u,v) \leq \|f(u)-f(v)\|_2 \leq d_G(u,v).
    \end{equation}
    This directly verifies the bi-Lipschitz hypothesis of Theorem~\ref{Thm_F_isTheFinalCoalgebra}. The completed coalgebra on $\Ss$ is given by 
    \begin{equation} 
    \nonumber
        \g_{\Ss}(z)=i\otimes(2z-\alpha_i) \quad \text{when } z\in\f_i(\Ss).
    \end{equation}
    Thus, the comparison theorem recovers the known finality in the set and continuous metric categories.
\end{example}

\begin{example}[$p$-adic integers $\mathbb Z_p$]
\label{Ex_p_adic}
    Final-coalgebra presentations of $\mathbb Z_p$ appear in \cite{Bhatpadic,Mano}, with \cite{Mano} establishing the completion construction and finality for continuous maps.
    
    For a prime $p$, put $V_p=\{0,\ldots,p-1\}$ for the digit alphabet and take $k=1$, copy set $M=[p]$, $D=((1,1))$, and $R$ the equality relation. Then
    \begin{equation} 
    \nonumber
        \F(X,x_0)=\bigl([p]\times X,(1,x_0)\bigr),
    \end{equation}
    with no gluing between copies. Metric admissibility is automatic in this one-pointed case.
    
    The set-theoretic initial algebra has carrier
    \begin{equation} 
    \nonumber
     G=\{s\in V_p^\omega:
           s_n=0\text{ for all sufficiently large }n\},
    \end{equation}
    pointed by $0^\omega$, with algebra map
    \begin{equation} 
    \nonumber
        g:[p]\times G\longrightarrow G,\qquad g(i,s)=(i-1)s,
    \end{equation}
    where copy $i\in[p]$ prepends the digit $i-1$. For $\rho=1/p$, its canonical metric is
    \begin{equation} 
    \nonumber
        d_G(s,t)=p^{-\ell(s,t)}\quad(s\ne t), \qquad d_G(s,s)=0,
    \end{equation}
    where $\ell(s,t)$ is the length of the longest common prefix.
    
    The $p$-adic evaluation defines a pointed bijection
    \begin{equation} 
    \nonumber
        f:G\longrightarrow\mathbb N_0,\qquad  f(s)=\sum_{r=0}^{\infty}s_rp^r,
    \end{equation}
    where the sum is finite because $s$ is eventually zero. Here $\mathbb N_0 = \{0,1,2,\cdots\}$. It satisfies
    \begin{equation} 
    \nonumber
        f(g(i,s))=(i-1)+pf(s) \quad\text{ and } \quad |f(s)-f(t)|_p=d_G(s,t).
    \end{equation}
    Thus it identifies the symbolic initial algebra isometrically with
    \begin{equation} 
    \nonumber
        [p]\times\mathbb N_0\longrightarrow\mathbb N_0, \qquad (i,n)\longmapsto (i-1)+pn,
    \end{equation}
    pointed by $0$ and equipped with the $p$-adic metric.
    
    The Cauchy completion of $G$ can be identified with $V_p^\omega$, carrying the same prefix metric and distinguished point $0^\omega$. The evaluation map extends to an isometric bijection
    \begin{equation} 
    \nonumber
        \widehat f:V_p^\omega \longrightarrow \mathbb Z_p, \qquad \widehat f(s)=\sum_{r=0}^{\infty}s_rp^r.
    \end{equation}
    The completed final coalgebra is the head--tail map
    \begin{equation} 
    \nonumber
        s_0s_1s_2 \cdots \longmapsto (s_0+1)\otimes(s_1s_2\cdots).
    \end{equation}
    Under $\widehat f$, this becomes the least-significant-digit coalgebra
    \begin{equation} 
    \nonumber
        \gamma(z)=(a_0+1)\otimes\frac{z-a_0}{p}, \qquad z\equiv a_0\pmod p,\quad  a_0\in V_p.
    \end{equation}
    Theorems~\ref{Thm_FinalCoalgebraInSet} and \ref{Thm_FinalCoalgebraInMet} therefore recover the known finality in the one-pointed set and continuous metric categories. 
 \end{example}

\begin{example}[$p$-adic integers $\mathbb Z_p$ : $p$-pointed IFS realization]
\label{Ex_p_adicIFS}
    We apply the IFS construction of Definition~\ref{Def_RelationForIFS}, distinguishing all branch fixed points. This gives a $p$-pointed presentation, differing from the one-pointed completion presentation in Example \ref{Ex_p_adic}.
    
    For a prime $p$, put $V_p=\{0,\ldots,p-1\}$ for the digit alphabet and consider the IFS
    \begin{equation} 
    \nonumber
        \f_i(z)=(i-1)+pz,\qquad i\in[p],
    \end{equation}
    on $\mathbb Z_p$. These maps satisfy
    \begin{equation} 
    \nonumber
        d_p(\f_i(z),\f_i(w))=p^{-1}d_p(z,w)
    \end{equation}
    and have fixed points $\alpha_i=(i-1)(1-p)^{-1}$, which are distinct and pairwise at distance one. We invite the reader to verify that the orbit is 
    \begin{equation} 
    \nonumber
        O=\left\{ N+\frac{bp^n}{1-p}\,:\, n\in\mathbb N_0,\; N\in\{0, \ldots,p^n-1\},\; b\in V_p \right\} = \frac{1}{p-1}\mathbb Z.
    \end{equation}
    
    The associated relation satisfies
    \begin{equation} 
    \nonumber
        (i,j)R(i',j') \quad \Longleftrightarrow \quad (i-1)+p\alpha_j=(i'-1)+p\alpha_{j'}. 
    \end{equation}
    Reduction modulo $p$ gives $i=i'$, and then $j=j'$. Thus $R$ is equality, so there is no gluing between copies. Consequently, the specification has $k=p$, copy set $M=[p]$, and diagonal output data $D=((1,1),\ldots,(p,p))$. Indexing distinguished points by $i\in[p]$, the associated functor is 
    \begin{equation} 
    \nonumber
        \F\bigl(A,(x_i)_{i\in[p]}\bigr) = \bigl([p]\times A,((i,x_i))_{i\in[p]}\bigr).
    \end{equation}
    The specification is metrically admissible for every $0<\rho<1$, since its distinguished points lie in distinct, unglued copies. 
    
    The set-theoretic initial algebra has as its carrier the eventually constant streams over $V_p$,
    \begin{equation} 
    \nonumber
        G=\{s\in V_p^\omega : s\text{ is eventually constant}\},
    \end{equation}
    with $i$-th distinguished point $(i-1)^\omega$ and algebra map $g(i,s)=(i-1)s$, where copy $i$ prepends the digit $i-1$. For $\rho=1/p$, its canonical metric is
    \begin{equation} 
    \nonumber
        d_G(s,t)=p^{-\ell(s,t)}\quad(s\ne t), \qquad d_G(s,s)=0,
    \end{equation}
    where $\ell(s,t)$ is the length of the longest common prefix.
    
    The comparison map is the $p$-adic evaluation:
    \begin{equation} 
    \nonumber
        f:G\longrightarrow\mathbb Z_p,\qquad  f(s)=\sum_{r=0}^{\infty}s_rp^r.
    \end{equation}
    It sends $(i-1)^\omega$ to $\alpha_i$ and satisfies $f(g(i,s))=\f_i(f(s))$. For a finite prefix followed by a constant tail, it gives
    \begin{equation} 
    \nonumber
        f(a_0\cdots a_{n-1}b^\omega)  =  \sum_{r=0}^{n-1}a_rp^r+\frac{bp^n}{1-p} \in O.
    \end{equation}
    
    Thus $f(G)=O$, while uniqueness of $p$-adic expansions makes $f$ injective. Hence $f:G\to O$ is bijective. If $s,t$ first differ at position $\ell$, then
    \begin{equation} 
    \nonumber
        f(s)-f(t) = p^\ell\bigl(s_\ell-t_\ell+p\xi\bigr), \qquad \xi\in\mathbb Z_p.
    \end{equation}
    Here, $s_\ell-t_\ell+p\xi$ is a $p$-adic unit, so
    \begin{equation} 
    \nonumber
        d_p(f(s),f(t))=p^{-\ell}=d_G(s,t).
    \end{equation}
    Hence the comparison is an isometry. The orbit algebra
    \begin{equation} 
    \nonumber
        \varphi:[p]\times O\longrightarrow O,\qquad \varphi(i,z)=(i-1)+pz,
    \end{equation}
    is injective, since the branches are injective and their images lie in distinct residue classes modulo $p$.
    
    Since $O$ contains $\mathbb N_0$, it is dense in $\mathbb Z_p$. Thus $f$ extends to an isometric bijection $C(G)\to\mathbb Z_p$, and Theorem~\ref{Thm_F_isTheFinalCoalgebra} gives finality in the set and continuous metric categories of non-degenerate $p$-pointed objects. The transported coalgebra reads off the least significant digit:
    \begin{equation} 
    \nonumber
        \gamma(z)=(a_0+1) \otimes \frac{z-a_0}{p} \in \Fww{1/p} \mathbb Z_p,      \qquad z\equiv a_0\pmod p,\quad a_0\in V_p.
    \end{equation}
    It preserves the distinguished points because $\gamma(\alpha_i)=i\otimes\alpha_i$.
\end{example}

\begin{example}[Nested radicals]\label{Ex_Nested_Radicals_IFS}
    Nested radicals of the form
    \begin{equation} 
    \nonumber
       \epsilon_1\sqrt{ 2+\epsilon_2\sqrt{ 2+\epsilon_3\sqrt{2+\cdots} }}, \qquad        \epsilon_i\in\{-1,1\},
    \end{equation}
    are governed by the same two-copy gluing specification as the unit interval. 
    
    Regard $X=[-2,2]$ as a two-pointed metric space with distinguished points $-2$ and $2$, equipped with the Euclidean metric scaled by a factor of $1/4$ so that it is one-bounded. Define
    \begin{equation} \nonumber
       \psi_{-1}(x)=-\sqrt{2-x}, \qquad \psi_1(x)=\sqrt{2+x}, \qquad x\in X.
    \end{equation}
    Equivalently, $ \psi_\delta(x)=\delta\sqrt{2+\delta x}, \qquad \delta\in\{-1,1\}$.     The maps $\psi_{-1}$ and $\psi_1$ take $X$ homeomorphically onto $[-2,0]$ and $[0,2]$, respectively. As the overlap agrees with the gluing relation of the two-copy gluing functor of the unit interval (cf. Example \ref{Ex_unitInterval}), they induce a well-defined continuous bijection
    \begin{equation} 
    \nonumber
       \psi:\Fww{\frac12}X\longrightarrow X, \qquad \psi(m\otimes x)=\psi_{2m-3}(x).
    \end{equation}
    Here $\Sigma$ has $M=[2]$, as in the unit interval, and the sign $2m-3$ corresponds to copy $m\in[2]$.
    
    Since $\Fww{\frac12}X$ is compact and $X$ is Hausdorff, $\psi$ is a homeomorphism. Its inverse is the continuous coalgebra structure
    \begin{equation} 
    \nonumber
       \psi^{-1}:X\longrightarrow\Fww{\frac12}X, \qquad
       \psi^{-1}(x)=
       \begin{cases}
          1\otimes(2-x^2), & -2\leq x\leq0,\\[1mm]
          2\otimes(x^2-2), & 0\leq x\leq2.
       \end{cases}
    \end{equation}
    in $\mcndk{2}$. In particular, even though the radical branches are not contractions, the continuous-coalgebra theorem is still applicable.
    
    Consider the balanced unit interval, equipped with the Euclidean metric scaled by $1/2$ to make it one-bounded, with 
    \begin{equation} 
    \nonumber
       \phi_\delta(u)=\frac{u}{2}+\frac{\delta}{2}, \qquad \delta\in\{-1,1\}.
    \end{equation}
    The inverse of the induced homeomorphism
    \begin{equation} 
    \nonumber
       \phi:\Fww{\frac12}[-1,1]\longrightarrow[-1,1],\qquad \phi(m\otimes u)=\phi_{2m-3}(u)
    \end{equation}
    is the final coalgebra structure (as in Example~\ref{Ex_unitInterval_IFS}). Therefore, there is a unique coalgebra morphism $ f:[-2,2]\longrightarrow[-1,1]$ satisfying
    \begin{equation} 
    \nonumber
       (\Fww{\frac12}f)\circ\psi^{-1} = \phi^{-1}\circ f.
    \end{equation}
    Define
    \begin{equation} 
    \nonumber
       h(x)=\frac{2}{\pi}\arcsin\left(\frac{x}{2}\right).
    \end{equation}
    The half-angle identities give
    \begin{equation} 
    \nonumber
       h\bigl(\psi_\delta(x)\bigr) = \phi_\delta\bigl(h(x)\bigr) = \frac{h(x)}{2}+\frac{\delta}{2}.
    \end{equation}
    Equivalently,
    \begin{equation} 
    \nonumber
       (\Fww{\frac12}h)\circ\psi^{-1} = \phi^{-1}\circ h.
    \end{equation}
    Thus $h$ is a coalgebra morphism, and finality implies $f=h$. Moreover, $h$ is a homeomorphism with inverse
    \begin{equation} 
    \nonumber
       h^{-1}(u)=2\sin\left(\frac{\pi u}{2}\right).
    \end{equation}
    Consequently, $([-2,2],\psi^{-1})$ is another realization of the final coalgebra.
    
    To evaluate the radicals, let $\delta_i=\epsilon_1\epsilon_2\cdots\epsilon_i$. For the finite radical 
    \begin{equation} 
    \nonumber
       x_r = \epsilon_1\sqrt{ 2+\epsilon_2\sqrt{ 2+\cdots+\epsilon_r\sqrt2 }},
    \end{equation}
    one has $x_r = \psi_{\delta_1} \bigl( \psi_{\delta_2} (\cdots\psi_{\delta_r}(0)\cdots) \bigr)$. The conjugacy therefore gives
    \begin{equation} 
    \nonumber
       h(x_r) = \phi_{\delta_1} \bigl( \phi_{\delta_2} (\cdots\phi_{\delta_r}(0)\cdots)        \bigr) = \sum_{i=1}^{r}\frac{\delta_i}{2^i}.
    \end{equation}
    
    Using the continuity of $h^{-1}$ and passing to the limit yields
    \begin{equation} 
    \nonumber
       \epsilon_1\sqrt{ 2+\epsilon_2\sqrt{ 2+\epsilon_3\sqrt{2+\cdots} } } =        2\sin\left( \frac{\pi}{2} \sum_{i=1}^{\infty} \frac{\epsilon_1\epsilon_2\cdots\epsilon_i}{2^i} \right).
    \end{equation}
    Conversely, every $u\in[-1,1]$ admits a signed binary expansion
    \begin{equation} 
    \nonumber
       u=\sum_{i=1}^{\infty}\frac{\delta_i}{2^i}, \qquad \delta_i\in\{-1,1\}.
    \end{equation}
    Setting
    \begin{equation} 
    \nonumber
       \epsilon_1=\delta_1, \qquad \epsilon_i=\delta_{i-1}\delta_i \quad (i\geq2)
    \end{equation}
    produces the corresponding nested radical. Thus the convergence and classical trigonometric evaluation of these radicals follow from the final coalgebra correspondence; compare \cite{nested}.
\end{example}
  
\begin{example}[The Sierpiński carpet: a limitation of finite \textit{point} gluing]
\label{Ex_Sierpinski_Carpet_IFS}
    
    The standard Sierpiński carpet provides a useful limitation of the finite
    ``{point}-gluing" framework, as it requires gluing along line segments and suggests open avenues for subsequent study. 
    
    Let $\Lambda=\{0,1,2\}^2\setminus\{(1,1)\}$, indexed in lexicographic order by
    \begin{equation}
    \nonumber
        (\lambda_1,\ldots,\lambda_8)=\bigl((0,0),(0,1),(0,2),(1,0),(1,2),(2,0),(2,1),(2,2)\bigr).
    \end{equation}
    For $i\in[8]$, define
    \begin{equation} 
    \nonumber
       \f_i(x)=\frac{x+\lambda_i}{3}, \qquad x\in[0,1]^2.
    \end{equation}
    Equip $[0,1]^2$ with the box metric, so that it is one-bounded. We have $\Delta = 1/2$ and $\lambda =1/3$.  
    
    The fixed point of $\f_i$ is $\alpha_i=\lambda_i/2$. Thus $O_0=\{\alpha_i:i\in[8]\}$ consists of eight distinct points and is therefore non-degenerate. Let $O$ be the orbit.  The attractor of this IFS is the Sierpiński carpet $\mathbb S_{\mathrm{car}} = \overline O $. 
     
    The associated gluing specification, with copy set $M=[8]$, is 
    \begin{equation}
    \nonumber
        \Sigma=\bigl([8],D=((1,1),\ldots,(8,8)),R\bigr),
    \end{equation}
    where  the relation $R$ records the coincidences among the distinguished points in the first-level pieces;
    \begin{equation} \nonumber
       (i,j)R(i',j')\quad\Longleftrightarrow\quad\f_i(\alpha_j)=\f_{i'}(\alpha_{j'}).
    \end{equation}
    Equivalently, $2\lambda_i+\lambda_j=2\lambda_{i'}+\lambda_{j'}$.
    
    Since $\lambda<\Delta$, the outputs are copy-separated. Choose $\rho\in[1/2,1)$; the specification is then metrically admissible. The associated orbit algebra is
    \begin{equation} 
    \nonumber
       \varphi:\Fw O\longrightarrow O, \qquad \varphi(i\otimes x)=\f_i(x).
    \end{equation}
    Finite-stage evaluation gives a short map $f:G\to O$, since     $1/3\leq\rho$. The estimate \eqref{Eqn_Mi_inequality} gives finite covers of $G$ by sets of arbitrarily small diameter, so $C(G)$ is compact. Consequently,
    \begin{equation} 
    \nonumber
       \ft  =\kappa_O\circ C(f):C(G)\longrightarrow\mathbb S_{\mathrm{car}}
    \end{equation}
    is surjective because its compact image contains the dense orbit $O$.

    Although $\varphi$ is surjective, it is not injective, as the explicit kernel pair below shows.  Indeed, adjacent first-level pieces meet along whole line segments. The relation $\overline R$, however, records only finitely many points on such intersections. The additional identifications needed are described precisely by the kernel equivalence relation
    \begin{equation} 
    \nonumber
       K_\varphi = \F O\times_O\F O = \left\{ (u,v)\in\F O\times\F O \;\middle|\;      \varphi(u)=\varphi(v) \right\}.
    \end{equation}
    In terms of representatives, for $i,j\in[8]$ and $x,y\in O$, we have $\bigl(i\otimes x,j\otimes y\bigr)\in K_\varphi$ if and only if $\f_i(x)=\f_j(y)$. Since $\f_i(x)=(x+\lambda_i)/3$, this can be written explicitly as
    \begin{equation} 
    \nonumber
       K_\varphi = \left\{ \bigl(i\otimes x,j\otimes y\bigr) \;\middle|\; i,j\in[8],\ x+\lambda_i=y+\lambda_j \right\}.
    \end{equation}
    For example, $\bigl(1\otimes p,4\otimes q\bigr)\in K_\varphi$, where $p=\f_6(\alpha_8)=(1,1/3)$ and $q=\f_1(\alpha_3)=(0,1/3)$. These are distinct quotient points because neither $p$ nor $q$ is a distinguished point of $O$. Since $\varphi$ is surjective, it induces a bijection
    \begin{equation} 
    \nonumber
       \overline\varphi: \frac{\F O}{K_\varphi} \xrightarrow{\ \cong\ } O, \qquad    \overline\varphi([i\otimes x])=\f_i(x).
    \end{equation}
    Thus the orbit can be recovered set-theoretically by quotienting by the full orbit-overlap relation. This quotient does not, however, identify $O$ with the initial algebra of the original finite gluing functor.  The relation $K_\varphi$ depends on the entire orbit and contains identifications arising at arbitrarily high stages; it is not generated by the prescribed finite boundary relation $\overline R$.  As the injectivity of the orbit algebra in IFS checklist fails, Proposition~\ref{pro_IsoO_o_and_G_inSet} and the Categorical Hutchinson theorem of Section~\ref{Sec:IFS} \textit{do not apply} to this specification.  The obstruction is already set-theoretic, so it cannot be removed by changing the metric parameter $\rho$.

    In \cite{Sierpinsky_carpet}, a functor on square metric spaces records the segment identifications and yields a final coalgebra bi-Lipschitz equivalent to the classical carpet. This suggests extending the present initial-algebra completion construction and quantitative IFS comparison to gluing along more general boundary spaces.
\end{example}

\begin{appendices}
\section{The module associated with a finite gluing specification}
\label{Subsec:RigorousModuleConstruction}

Let $\mathcal A_k$ be the category with objects $0,1$ and precisely $k$ non-identity arrows $\epsilon_j:0\to1$, $j\in[k]$. Write $\mathcal P_k$ for the full subcategory of $[\mathcal A_k,\mathbf{Set}]$ whose objects have singleton $0$-component. The equivalence $E:\set\simeq\mathcal P_k$ is given by
\begin{equation} 
\nonumber
    EA(0)=\{*\},\qquad EA(1)=A,\qquad EA(\epsilon_j)(*)=a_j,
\end{equation}
with $E(f)_0=\operatorname{id}_{\{*\}}$ and $E(f)_1=f$. Conversely, $X\in\mathcal P_k$ determines the pointed set $(X(1),X(\epsilon_1)(*),\ldots,X(\epsilon_k)(*))$.

Fix a finite gluing specification $\Sigma=(M,D,R)$, with $M=[\s]$ and output labels $d_r=(i_r,j_r)$. Let $\mathcal C_\Sigma$ denote the underlying set of the pointed set $\mathcal C=\F I$. Thus
\begin{equation}
\nonumber
    \mathcal C_\Sigma=(M\times I)/\overline R,
    \qquad i\otimes j=[(i,j)]_{\overline R},\qquad C_r=\eta(r)=i_r\otimes j_r.
\end{equation}

\begin{proposition}[The associated module]
\label{Prop:RigorousAssociatedModule}
    The specification $\Sigma$ determines a finite module
    \begin{equation} 
    \nonumber
        \mathcal M_\Sigma:\mathcal A_k^{op}\times\mathcal A_k      \longrightarrow\mathbf{Set}
    \end{equation}
    whose tensor functor restricts to $\mathcal P_k$ and satisfies the natural isomorphism
    \begin{equation} 
    \nonumber
        \mathcal M_\Sigma\otimes EA\cong E(\F A).
    \end{equation}
\end{proposition}

\begin{proof}
Define the four values by
    \begin{equation} 
    \nonumber
    \begin{array}{c|cc}
        \mathcal M_\Sigma(p,q)&q=0&q=1\\ \hline
        p=0&\{*\}&\mathcal C_\Sigma\\
        p=1&\emptyset&M.
    \end{array}
    \end{equation}
The non-identity actions with nonempty domain are
\begin{equation} 
\nonumber
\begin{aligned}
    \mathcal M_\Sigma(\epsilon_j,1)(i) & =i\otimes j,\\
    \mathcal M_\Sigma(0,\epsilon_r)(*)&=C_r=\eta(r).
\end{aligned}
\end{equation}
All remaining non-identity actions have empty domain and are unique. There are no composites of non-identity arrows in $\mathcal A_k$, and every square comparing the two variable actions has domain $\mathcal M_\Sigma(1,0)=\emptyset$. Thus these data define a profunctor. It is finite because $\mathcal A_k$ and all four values are finite.

For a pointed set $A$, the coend at $0$ is a singleton, while at $1$ it is
\begin{equation} 
\nonumber
    (\mathcal M_\Sigma\otimes EA)(1) \cong \frac{\mathcal C_\Sigma\amalg(M\times A)} {i\otimes j\sim(i,a_j)}
    \cong\F A.
\end{equation}
Indeed, eliminating the $\mathcal C_\Sigma$ summand imposes exactly the boundary identifications generated by $R$. The arrow $\epsilon_r$ selects $C_r$, which corresponds to $i_r\otimes a_{j_r}$, the prescribed $r$-th distinguished point. For a pointed map $f:A\to B$, the induced carrier map sends $i\otimes a$ to $i\otimes f(a)$. Hence the identification preserves the pointed structure and is natural in $A$.
\end{proof}

\begin{proposition}[Module non-degeneracy criterion]
\label{Prop_RNondegenerate_ModuleNondegenerate}
    The module $\mathcal M_\Sigma$ is non-degenerate if and only if $\Sigma$ satisfies \eqref{eqn:Non-degeneracy}. These conditions are therefore equivalent to the restriction criterion for $\F-$ on $\setnd$ in Proposition~\ref{Prop_EndofunctorOnSetND}.
\end{proposition}

\begin{proof}
By \cite[Definition~2.7]{Lein}, module non-degeneracy means non-degeneracy of each row $\mathcal M_\Sigma(p,-)$. For $\mathcal A_k$, the criterion of \cite[Theorem~4.11]{Lein} says that a functor $X:\mathcal A_k\to\mathbf{Set}$ is non-degenerate precisely when the maps $X(\epsilon_j)$ are injective with pairwise disjoint images. Here ND1 applied to equal or distinct parallel arrows gives these two requirements, and ND2 adds no further condition; compare \cite[Example~4.2]{Lein}.

The row $\mathcal M_\Sigma(1,-)$ has empty $0$-component, so the criterion holds automatically. In the row $\mathcal M_\Sigma(0,-)$, every boundary map has singleton domain and image $\{C_r\}$. Thus this row is non-degenerate exactly when the classes $C_1,\ldots,C_k$ are pairwise distinct, which is \eqref{eqn:Non-degeneracy}. The final assertion follows from Proposition~\ref{Prop_EndofunctorOnSetND}.
\end{proof}

\section{Preservation of completeness}
\label{App_FPreservesCompleteness}

Here, we shall prove Lemma \ref{Lem_FPreservesCompleteness}, that for non-degenerate $\Sigma$, $\Fw -$ preserves complete metric spaces. 

\begin{proof}[Proof of Lemma \ref{Lem_FPreservesCompleteness}]
Non-degeneracy follows from Proposition \ref{Prop_EndofunctorOnSetND}. We need to show the convergence of Cauchy sequences.

Let $(z_n)_{n\geq1}$ be a Cauchy sequence in ${\Fw A}$.  We distinguish two cases.

\smallskip
\noindent\emph{Case 1: $\liminf_{n\to\infty}r(z_n)=0$.}
Choose a subsequence $(z_{n_j})$ such that
\begin{equation} 
\nonumber
   r(z_{n_j})<\frac1j.
\end{equation}
For every $j$, choose $p_j\in \mathcal C_A$ satisfying $d_{\Fw A}(z_{n_j},p_j)=r(z_{n_j})$.  Since $\mathcal C_A$ is finite, some $p\in \mathcal C_A$ occurs along a further subsequence.  After passing to and relabelling that subsequence, we therefore have
\begin{equation} 
\nonumber
   \dqm{A}(z_{n_j},p)=r(z_{n_j})\longrightarrow0.
\end{equation}
Thus $(z_n)$ has a subsequence converging to $p$.  Because $(z_n)$ is Cauchy, the entire sequence converges to $p$.  

\smallskip
\noindent\emph{Case 2: $\liminf_{n\to\infty}r(z_n)>0$.}
There are $\varepsilon_0>0$ and $N_0$ such that
\begin{equation}
   r(z_n)\geq\varepsilon_0    \qquad \text{ for } n\geq N_0.
\label{eqn:bound_4}
\end{equation}
Put $\varepsilon=\min\{1,2\varepsilon_0\}>0$. Since $(z_n)$ is Cauchy, after increasing $N_0$, if necessary, we may also assume that
\begin{equation}  
    \dqm{A}(z_n,z_\ell)<\varepsilon \qquad(n,\ell\geq N_0).
\label{eqn:bound_5}
\end{equation}

In particular, $z_n\notin \mathcal C_A$ for $n\geq N_0$, so write its unique representative as $ z_n=m_n\otimes x_n$, with $x_n\notin A_{\partial}$. If $m_n\neq m_\ell$ for some $n,\ell\geq N_0$, then Lemma \ref{lem:metric_behaviour_away_form_guling_locus}-\eqref{eqn:bound_1} and \eqref{eqn:bound_4} give
\begin{equation} \nonumber
   \dqm{A}(z_n,z_\ell)  \geq\min\{1,r(z_n)+r(z_\ell)\} \geq\varepsilon,
\end{equation}
contrary to \eqref{eqn:bound_5}.  Hence, all terms in this tail lie in one fixed copy. Say $z_n=m\otimes x_n$ for $n\geq N_0$. Now, for $n,\ell\geq N_0$, equations \eqref{eqn:bound_4} and \eqref{eqn:bound_5} imply
\begin{equation} 
\nonumber
   \dqm{A}(z_n,z_\ell) <\varepsilon \leq\min\{1,r(z_n)+r(z_\ell)\}.
\end{equation}
Applying Lemma \ref{lem:metric_behaviour_away_form_guling_locus}-\eqref{eqn:bound_2} yields
\begin{equation} \nonumber
   d_A(x_n,x_\ell)=\rho^{-1}d_{\Fw A}(z_n,z_\ell).
\end{equation}
It follows that $(x_n)_{n\geq N_0}$ is Cauchy in $A$.  Since $A$ is complete, there is $x\in A$ such that $x_n\to x$.  A one-step chain in the $m$-copy gives
\begin{equation} \nonumber
   \dqm{A}\bigl(m\otimes x_n,m\otimes x\bigr) \leq\rho \cdot d_A(x_n,x)\longrightarrow0.
\end{equation}
Therefore $z_n\to m\otimes x$ in ${\Fw A}$.

In both cases the Cauchy sequence $(z_n)$ converges in $\Fw A$.  Hence
$M\otimes_\rho A$ is complete.
\end{proof}
\end{appendices}

\bibliographystyle{alpha}
\bibliography{sn-bibliography}

\end{document}